\documentclass[a4]{amsart}

\usepackage{latexsym}
\usepackage{verbatim}

\newcommand{\R}{{\mathbb R}}

\newcommand{\B}{{\mathbb B}}

\newcommand{\half}{{\frac{1}{2}}}

\newcommand\Ha{\operatorname{Ha}^{(1)}}
\newcommand\Id{\operatorname{Id}}
\newcommand\cTY{\tilde T^* Y}
\newcommand\WF{\operatorname{WF}'}
\newcommand\csd{\Delta_{S^* Y}}
\newcommand\Cbill{C_{\operatorname{billiard}}}

\newcommand\Omegabar{\overline{\Omega}}

\newtheorem{theo}{{\sc Theorem}}[section]
\newtheorem{cor}[theo]{{\sc Corollary}}
\numberwithin{equation}{section}

\newtheorem{lem}[theo]{{\sc Lemma}}
\newtheorem{prop}[theo]{{\sc Proposition}}

\newtheorem{defn}[theo]{{\sc Definition}}

\newenvironment{rem}{\medskip\noindent{\it Remark.\/} }{\medskip}

\newcommand\pa{\partial}

\newcommand\vol{\operatorname{vol}}
\newcommand\tr{\operatorname{tr}}
\newcommand{\RR}{{\mathbb R}}
\newcommand\Op{\operatorname{Op}}
\newcommand\supp{\operatorname{supp}}
\def\dbyd#1#2{\frac{ \partial #1}{\partial #2}}
\newcommand\Neu{{\operatorname{Neu}}}
\newcommand\Dir{\operatorname{Dir}}

\newcommand\llangle{\langle \langle}
\newcommand\rrangle{\rangle \rangle}

\newcommand\pt{\vskip 5pt \noindent$\bullet$ \ }

\newcommand\norm{\operatorname{norm}}

\title[Quantum ergodicity of boundary values of eigenfunctions]
{Quantum ergodicity of boundary values of eigenfunctions}

\author{Andrew Hassell}
\address{Department of Mathematics, Australian National University, Canberra 0200 ACT Australia}
\email{hassell@maths.anu.edu.au}
\author{Steve Zelditch}
\address{Department of Mathematics, Johns Hopkins University, Baltimore,
MD
21218, USA}
\email{ zelditch@math.jhu.edu}

\thanks{The first author was partially supported by an Australian Research Council Fellowship. The second author was partially supported by NSF grant  \#DMS-0071358.}

\begin{document}

\maketitle


\begin{abstract} Suppose that $\Omega \subset \RR^n$ is a bounded, piecewise smooth domain.
  We prove that the boundary values (Cauchy data)
of eigenfunctions of the Laplacian on $\Omega$ with various
boundary conditions are quantum ergodic if the classical billiard
map $\beta$ on the ball bundle $B^*(\partial \Omega)$ is ergodic.

 Our proof is based on the classical
observation that the boundary values of an interior eigenfunction
$\phi_{\lambda}$, $\Delta \phi_{\lambda} = \lambda^2 \phi_\lambda$
is an eigenfunction of an operator $F_{h}$ on the boundary of
$\Omega$ with $h = \lambda^{-1}$. In the case of the Neumann
boundary condition, $F_h$ is the boundary integral operator
induced by the  double layer potential. We show that $F_{h}$ is a semiclassical Fourier integral operator
quantizing the billiard map plus a `small' remainder; the quantum
dyanmics defined by $F_h$ can be exploited on the boundary  much
as the quantum dynamics generated by the wave group were exploited
in the interior of domains with corners and ergodic billiards in
the work of Zelditch-Zworski (1996). Novelties include the facts
that $F_h$ is not unitary and (consequently) the boundary values
are equidistributed by measures which are not invariant under
$\beta$ and which depend on the boundary conditions.

Ergodicity of boundary values of eigenfunctions on domains with
ergodic billiards  was conjectured by S. Ozawa (1993), and was
almost simultaneously proved by Gerard-Leichtnam (1993) in the
case of convex $C^{1,1}$ domains (with continuous tangent planes)
and with Dirichlet boundary conditions. Our methods seem to be
quite different. Motivation to study piecewise smooth domains comes from the fact that almost all known ergodic domains are of this form.

\end{abstract}

\tableofcontents

\section{Introduction}

The purpose of this article is to prove quantum ergodicity
\begin{equation} \langle A_{h_j} u_j^b, u_j^b \rangle  \to \int_{B^* Y} \sigma(A) \; d
\mu_B , \ j \to \infty \  \mbox{along a density one sequence} \;\;
\label{q-erg}\end{equation}
of the  boundary values $u_j^{b}$  of interior eigenfunctions
$$\left\{ \begin{array}{l} \Delta \; u_j = \lambda_j^2\;  u_j\;\; \mbox{in} \; \Omega, \;\; \langle u_j, u_k
 \rangle_{L^2(
\Omega)} = \delta_{jk}, \\ \\B u_j |_{Y} = 0,\;\;\;\; Y = \partial
\Omega
\end{array} \right.
$$ of the Euclidean  Laplacian $\Delta$ on a compact piecewise smooth domain $\Omega \subset
\R^n$ and with  classically ergodic
billiard map $\beta:  B^* Y \to B^*Y$, where $Y = \partial \Omega$.
Here $A_h$ is a zeroth order semiclassical pseudodifferential operator on $Y$.
The relevant notion of boundary values  (i.e. Cauchy data)
$u_j^{b}$  depends on the boundary condition $B$, as does the
classical limit measure $d\mu_B$ according to which the boundary
values become equidistributed. Our methods cover  Dirichlet,
Neumann, Robin and more general boundary conditions of the
form
\begin{equation} B u = \partial_{\nu} u - K ( u|_{Y} ),\;\; K
\in \Psi^1(Y) \text{ self-adjoint with non-negative symbol } \label{K-bc}\end{equation}
which we refer to as a $\Psi^1$-Robin boundary condition.

For the interior problem,  ergodicity of eigenfunctions of
Laplacians on bounded domains with corners and  with ergodic
billiard flow was proved by Zelditch-Zworski \cite{ZZw},
following an earlier proof by Gerard-Leichtnam \cite{GL} in the
case of $C^{1,1}$  convex domains. Our proof of boundary
ergodicity is independent of these proofs in the interior case. In
the case of manifolds without boundary, results on ergodicity of
eigenfunctions originate in the work of A. I. Schnirelman \cite{Schnirelman} and were
carried forward by Zelditch \cite{Z8} and Colin de Verdiere \cite{CdeV}. We refer to
\cite{Z} for background and a simple proof which will be developed
here.

Ergodicity of boundary values of eigenfunctions was conjectured by
S. Ozawa \cite{Ozawa} in 1993 and was independently stated and
proved by Gerard-Leichtnam \cite{GL} in the same year in the case
of Dirichlet boundary conditions on $C^{1,1} $ convex domains
(i.e. domains whose unit normal is Lipschitz regular). Our
extension to nonconvex piecewise smooth domains is motivated by the fact
that `most' known ergodic billiard tables are nonsmooth. More
precisely, $C^2$ convex domains never have ergodic billiards, and
the only known
 $C^{1,1}$ convex ergodic billiard tables appear to be the  Bunimovich stadium example \cite{BU2, BU3}
 and its relatives  in two
dimensions (i.e. plane domains formed by straight segments and
circular segments).  In particular, no $C^{1,1} $ convex ergodic
domains are known in dimensions $\geq 3$ (the example cited in
\cite{GL} from \cite{BU3} is not ergodic) and there is some doubt
that they exist \cite{BU1}. The known higher dimensional ergodic stadia (i.e. with only
focussing or nowhere dispersing boundary faces, consisting of convex or flat components)  are non-convex
and have corners  \cite{BR,BR2}.
A sizable collection of ergodic,
piecewise smooth plane domains is furnished by generic polygons
\cite{KMS}. Many other  examples of non-convex ergodic (and hyperbolic) billiard domains in higher dimensions are given by
dispersing billiard domains bounded by unions of concave boundary components (see e.g.
\cite{W}).

Our proof is based on a reduction to the boundary of the
eigenvalue problem.  The intuitive idea is that the Cauchy data
$(u_j|_Y, \partial_{\nu} u_j |_Y)$ of interior eigenfunctions
$u_j$ provide a kind of quantum cross section to the interior
eigenfunctions, just as the billiard map $\beta$ on $B^*Y$
provides a kind of cross section to the billiard flow on $T^*
\Omega.$  To be precise, our starting point is the
 classical observation that the boundary value of an interior
eigenfunction with eigenvalue $\lambda^2$ is an eigenfunction of a certain boundary operator
$F_{h}$ with $h = \lambda^{-1}$. We  analyse $F_h$ (on a convex
domain, or a modification of $F_h$ for a nonconvex domain) as a semiclassical Fourier integral operator quantizing
$\beta$ plus a remainder which is almost ignorable. Boundary
ergodicity has some new features which are not present in interior
ergodicity, stemming from the fact that the dynamics generated by
$F_h$ defines an  endomorphism but not an automorphism of the
observable algebra.

To state our results, we will need some notation:  Let $\Delta_B$
 denote the positive Laplacian on $\Omega$ with boundary
conditions $B u = 0$. Then $\Delta_B$ has discrete spectrum $0 <
\lambda_1 < \lambda_2 \leq \dots \to \infty$, where we repeat each
eigenvalue according to its multiplicity, and for each $\lambda_j$
we may choose an $L^2$ normalized eigenfunction $u_j$.

To each boundary condition $B$ corresponds
\begin{itemize}

\item  A specific notion of  boundary value $u_j^{b}$ of the eigenfunctions $u_j$.
We denote the $L^2$-normalized boundary values by $\hat u_j^{b} =
u_j^b/||u_j^b||.$

\item A specific measure $d\mu_B$ on $B^* (Y). $

\item A specific state $\omega_B$ on the space $\Psi_h^0(Y)$
of semiclassical pseuodifferential operators of order zero.
\end{itemize}

The correspondence is dictated by the local Weyl law (Lemma~\ref{LWL}) for the boundary condition $B$. Here
is a table of the relevant boundary value notions.  In the table,
$\kappa$ denotes a $C^{\infty}$ function on $Y$ while  $k$ is the
principal symbol of the operator $K \in \Psi^1(Y)$ in
\eqref{K-bc}, and $d \sigma$ is the natural symplectic volume
measure on $B^*Y$. We also define the function $\gamma(q)$ on
$B^*Y$ by
\begin{equation}
\gamma(q) = \sqrt{1 - |\eta|^2}, \quad q = (y,\eta).
\label{a-defn}\end{equation}

\bigskip
\noindent
\hskip 50pt
{\large \begin{tabular}{|c|c|c|c|c|} \hline
\multicolumn{4}{|c|}{\bf Boundary Values} \\ \hline
 B & $Bu$  &   $u^{b}$ & $d\mu_B$   \\
\hline
 Dirichlet &  $u|_{Y}$ &  $\partial_{\nu} u |_{Y}$ & $\gamma(q) d\sigma $ \\
\hline
 Neumann &   $\partial_{\nu} u |_{Y}$  &  $u|_{Y}$  & $\gamma(q)^{-1} d\sigma$  \\
\hline
 Robin  &  $(\partial_{\nu} u - \kappa u) |_{Y}$ &  $u|_{Y}$ & $\gamma(q)^{-1} d\sigma$  \\ \hline
$\Psi^1$-Robin & $(\partial_{\nu} u - K u) |_{Y}$ & $u|_{Y}$ & $\displaystyle{\frac{\gamma(q)d\sigma }{ \gamma(q)^2 + k(q)^2}} $\\
\hline
 \end{tabular}}
\bigskip

Throughout, we assume our domain $\Omega$ is a
piecewise smooth manifold embedded in $\R^n$. Hence its boundary,
denoted $Y$, is the union of a smooth part $Y^o$ and a
singular set $\Sigma$, which has measure zero (see section~\ref{corners} for
background on such manifolds). The metric on $\Omega$ is understood to be the
Euclidean metric.  Our main result is that, if the
billiard ball map $\beta$ on $B^*Y$ is ergodic, then the boundary
values $u_j^{b}$ of eigenfunctions are quantum ergodic. As
reviewed in \S 2, quantum ergodicity has to do with time and space
averages of observables. The relevant algebra of observables in
our setting is the algebra $\Psi_h^0(Y)$ of zeroth order
semiclassical pseudodifferential operators on $Y,$ depending on
the parameter $h \in [0, h_0]$. We denote the symbol of  $A = A_h
\in \Psi_h^0(Y)$ by  $a = a(y,\eta,h)$. Thus $a(y, \eta) = a(y,
\eta, 0)$ is a smooth function on $T^*Y$.
 We further define  states on the algebra $\Psi_h^0(Y)$ by
\begin{equation}\begin{aligned}
\omega_B(A) &= \frac{4}{\vol(S^{n-1})\vol(\Omega)} \int_{B^*Y} a(y,\eta)  d\mu_B.
\end{aligned}\label{omegaB}\end{equation}

Then our main result is
\begin{theo}\label{main}
Let  $\Omega \subset \RR^n$ be a bounded piecewise smooth manifold (see Definition~\ref{cpm}) with  ergodic billiard map. Let $\{u_j^{b}\}$
be the boundary values of the eigenfunctions
$\{u_j\}$ of $\Delta_B$ on $L^2(\Omega)$ in the sense of the table
above.
 Let $A_h$ be a semiclassical operator of order zero on $Y$. Then there is a subset $S$ of the positive integers, of density one,
such that
\begin{equation}\begin{gathered}
\lim_{j \to \infty, j \in S} \langle A_{h_j} u_j^b,
u_j^b \rangle = \omega_B(A), \quad B = \text{ Neumann, Robin or $\Psi^1$-Robin}, \\
\lim_{j \to \infty, j \in S} \lambda_j^{-2} \langle A_{h_j} u_j^b,
u_j^b \rangle = \omega_{B}(A), \quad B = \text{ Dirichlet},\end{gathered}
\label{main-eqn}\end{equation}
where $h_j = \lambda_j^{-1}$ and $\omega_B$ is as in \eqref{omegaB}. \end{theo}

\pt Let us give the results more explicitly for the identity operator. For the Neumann boundary condition,
$$
\lim_{j \to \infty, j \in S} \| u_j^b \|_{L^2(Y)}^2 = \frac{2 \vol(Y)}{\vol(\Omega)},
$$
while for the Dirichlet boundary condition,
$$
\lim_{j \to \infty, j \in S} \lambda_j^{-2} \| u_j^b \|_{L^2(Y)}^2 = \frac{2 \vol(Y)}{n\vol(\Omega)}.
$$

\pt As mentioned above, the case of Dirichlet boundary conditions
on $C^{1,1}$  convex domains was proved earlier by
G\'erard-Leichtnam \cite{GL}. Their proof was based on an identity
(\cite{GL}, Theorem 2.3) relating quantum limit measures of
interior Dirichlet eigenfunctions to those of its boundary values.
As they point out, the proof assumes  $C^{1,1}$ regularity and
does not apply to  domains with corners. Subsequent to the initial version of this
article (which only proved ergodicity for convex billiard domains),  N. Burq \cite{Burq} proved boundary
ergodicity for all Riemannian manifolds with corners (including non-convex domains) and
 all boundary conditions considered here. His method is to extend the method of \cite{GL} to
 general piecewise smooth domains and boundary conditions, and
 thus to reduce
 the  proof of boundary quantum ergodicity to the known interior quantum ergodicity result in \cite{ZZw}.
 His proof also uses the results of the present article on non-concentration of eigenfunctions at the
 corners.

\pt The Neumann and Dirichlet limit measures $d\mu_B$ can be
understood as follows: First,  $\gamma^{-1} d \sigma$ is the
projection to $B^* \Omega$ of the Liouville measure on the set
$S^*_{in} Y $ of inward pointing unit vectors to $\Omega$ along
$Y$ under the projection taking a vector to its tangential
component. We also recall that  $\beta$ is symplectic with respect
to the canonical symplectic form $d \sigma $  on $B^* Y$.  Thus,
boundary values of Neumann eigenfunctions are equidistributed
according to the measure on $B^* Y$ induced by interior Liouville
measure rather than the boundary symplectic volume measure.
In the case of Dirichlet boundary conditions, the boundary value $u_j^b$ is taken to be normal derivative of the eigenfunction at the boundary.
The symbol of $h \pa_\nu$, restricted to the spherical normal bundle, and then projected to $B^*Y$
is equal to $\gamma$, so we should expect to get the square of this factor in the Dirichlet case (since \eqref{q-erg} is quadratic in $u_j$) compared to the Neumann case. The normal derivatives in the Dirichlet case also account for the factor $\lambda_j^{-2}$ in \eqref{main-eqn}.

\vskip 5pt

The fact that the quantum limit state $\omega_B$ and the
corresponding measure $d\mu_B$ do not in general coincide with the
natural symplectic volume measure $d \sigma$ on $B^*Y$ will be
traced in  \S \ref{ENDO} to the fact that the quantum dynamics is
defined by an endomorphism rather than an automorphism of the
observable algebra. Let us explain how this works  in the case of
Neumann boundary conditions. In this case,  the dynamics are
generated by the operator $F_h$ on $Y$ with kernel
\begin{equation} \begin{gathered}
F_h(y,y') = 2\frac{\pa}{\pa \nu_y} G_0(y,y',h^{-1}), \quad y \neq y'
\in Y, \text{ where }\\
G_0(y,y',\lambda) =  \frac{i}{4} \lambda^{n-2} (2 \pi \lambda |z-z'|)^{-(n-2)/2}
\Ha_{n/2-1}(\lambda | z-z' |)
\end{gathered}\label{Neumann-F}\end{equation}
is the free outgoing Green function on $\RR^n$. By virtue of
Green's formula
\begin{equation} u_j(z) = \int_{}
\big[\partial_{\nu_{y'}} G_0(z, y', \lambda_j)  u_j(y')- G_0(z, y',\lambda_j) \partial_{\nu_{y'}} u_j(y')
 \big] d\sigma(y')
\end{equation}
for any solution of $\Delta u_j = \lambda_j^2 u_j,$  and the jump formula
\begin{equation}
\lim_{z \to y \in Y} \int_{Y} 2 \partial_{\nu_{y'}} G_0(z, y',
\lambda_j)  u_j(y') d\sigma(y') = u_j(y)  + F_h(u_j)(y),
\label{jumpformula}
\end{equation}
this operator
leaves the boundary values of Neumann eigenfunctions $u_j^{b}$
invariant:
\begin{equation} \label{FONE}
F_{h_j} u_j^{b} = u_j^{b}, \quad j = 1,2, \dots
\end{equation}
It follows that the states
\begin{equation} \rho_j(A) :=  \langle A_{h_j} u_j^b, u_j^b \rangle
\end{equation}
are invariant for $F_{h_j}$.  Similar invariance properties hold
for the other boundary conditions.  As we will show, the family
$\{F_{h}\}$ defines a semiclassical Fourier integral operator
associated to the billiard map $\beta$ (for convex $\Omega$), plus some terms which turn
out to be negligible for our problem. The quantum dynamics on
$\Psi^0_h (Y)$ is thus generated by the conjugation
\begin{equation} \label{QDF} \alpha_{h_j} (A_{h_j}) =
F_{h_j}^*\;A_{h_j} \; F_{h_j}
\end{equation}
This is analogous to the interior dynamics generated by
\begin{equation} \label{QD}
\alpha_t(A) = U_t A U_t^*,\;\;\; U(t) = e^{i t \Delta_B},
\end{equation} but it has one important difference: unlike $U(t)$,
$F_h$ is not unitary or even normal. Indeed, the zeroth order part of $F_h^* F_h$ is a
semiclassical pseudodifferential operator with a non-constant
symbol. This new feature of the quantum ergodicity problem is one
of the prinicipal themes of the present article.

We now outline the proof, emphasizing the aspects which are new to
the boundary case. The general strategy is the same as in \cite{Z,
ZZw}, and relies on two main ingredients: a local Weyl law for the
$u_j^b$, and an Egorov type theorem for an
  operator $F_h$. Naturally, the reduction to the boundary brings
  in additional considerations, which are of some independent
  interest.

We begin with the local Weyl law, which has nothing to do with
ergodicity; it is valid for all domains $\Omega$.

\begin{lem} \label{LWL}
 Let $A_h$ be either the identity operator on $Y$ or a zeroth order semiclassical operator on $Y$ with kernel supported away from the singular set.   Then for any of the above boundary conditions $B$,
 we have:
 \begin{equation}\begin{gathered}
\lim_{\lambda \to \infty} \frac{1}{N(\lambda)} \sum_{\lambda_j
\leq \lambda}  \langle A_{h_j} u_j^b, u_j^b \rangle =
\omega_B(A) , \quad  B = \text{ Neumann, Robin or $\Psi^1$-Robin}, \\
\lim_{\lambda \to \infty} \frac{1}{N(\lambda)} \sum_{\lambda_j
\leq \lambda} \lambda_j^{-2} \langle A_{h_j} u_j^b,
u_j^b \rangle = \omega_{B}(A), \quad B = \text{ Dirichlet}.\end{gathered}
\label{Weyl-limit}\end{equation}
\end{lem}

When $A$ is a multiplication operator, and for Dirichlet boundary
conditions, this local Weyl law was essentially proved by Ozawa
\cite{Ozawa}.  In \S 3 we extend the proof to general
semiclassical pseudodifferential operators and to the boundary
conditions described above. That allows us to capture uniform
distribution of eigenfunctions in phase space rather than just in
configuration space. For multiplication operators, the
local Weyl law can be obtained from the interior local Weyl law by
Hadamard's variational formula with respect to the boundary
conditions. This observation was first made by Ozawa \cite{O2, O3,
O4}.

From the local Weyl law we deduce an invariance property of the
limit states. For notational simplicity, we confine ourselves here to
Neumann boundary conditions, where   the boundary operator is
\eqref{Neumann-F}; analogous invariance properties hold for other boundary
conditions with small modifications to $F_h$.

\begin{cor} \label{INV} The state $\omega_{\Neu}$ is invariant under $F_h$: $\omega_{\Neu}(F_h^* A F_h) =
\omega_{\Neu}(A)$. \end{cor}

Indeed, the states $\rho_j(A) = \langle A_{h_j}
u_j^b, u_j^b \rangle$ are invariant so any average or
limit of averages of these states  will be invariant.

The Egorov type result for the operator $F_h$ is as follows:

\begin{lem} \label{Egorov} Let $\Omega \subset \RR^n$ be a bounded piecewise smooth convex domain, let $\beta$ denote the billiard map on
$B^* Y^o$ and let $A_h = \Op(a_h)$ be a zeroth order operator whose symbol $a(y,\eta,0)$ at $h=0$ is supported away from
$$ \{ |\eta| = 1 \} \cup \Sigma \cup \beta^{-1}(\Sigma).
$$
Let $\gamma$ be given by \eqref{a-defn}. Then
$$
F_h^* A_h  F_h = \tilde A_h + S_h,
$$
where $\tilde A_h$ is a zeroth order pseudodifferential operator
and $\| S_h \|_{L^2 \to L^2} \leq C h$. The symbol of $\tilde A_h$
is
\begin{equation}
 \tilde a = \begin{cases}
\gamma(q) \gamma(\beta(q))^{-1} a(\beta(q)) , \quad q \in B^*Y \\
0 , \phantom{\gamma(q)^{-1} \gamma(\beta(q)) b(\beta(q)) , \ } q
\notin B^*Y. \end{cases} \label{Egorov-formula}\end{equation}
\end{lem}

This is a rigorous version of the statement that $F_h$ quantizes
the billiard ball map. The unusual transformation law of the
symbol reflects the fact that \eqref{QDF} is not an
automorphism. This Egorov theorem is relevant to the Neumann
boundary problem. In the Dirichlet case, the relevant operator is
$F_h^*$. In the Robin case there is a lower order term, while in
the $\Psi^1$-Robin case there is a second term of the same order.
For nonconvex domains, we replace $F_h$ by a modified invariant operator which removes the spurious wavefront set of $F_h$ not associated with $\beta$.

We now sketch the completion of the proof in the case of  the Neumann
boundary condition. As in the case of automorphisms, it is
essentially a convexity argument (see \S 2 or \cite{Z} for this
point of view).  For simplicity of exposition, we  temporarily ignore the problems caused by the corners and pretend that
 the domain is smooth; details on the corner issues  appear in Section~\ref{Proof}.

To show that
$$
\langle A u_j^b, u_j^b \rangle \to \omega_{\Neu}(A),
$$
along a density one subsequence of integers $j$ is essentially to
show that
\begin{equation}\label{DESIRED}
\limsup_{\lambda \to \infty} \frac1{N(\lambda)} \sum_{\lambda_j <
\lambda}  \Big| \langle (A -  \omega_{\Neu}(A)) u_j^b,
u_j^b \rangle \Big|^2 = 0. \end{equation} Due to the novel
form of the local Weyl law and the Egorov theorem, we cannot apply
the averaging argument of \cite{Z} directly to (\ref{DESIRED}).
Instead, we  introduce an auxiliary pseudodifferential operator
$R$ which is almost invariant under conjugation by $F_h$. It
suffices that its symbol is given (approximately) by
\begin{equation}  \sigma_R(q) \sim c \gamma(q) = c (1 - |\eta|^2)^{1/2},\;\;\; c := \frac{\vol(S^{n-1}) \vol(\Omega)}{4 \vol(B^*Y)}.
\label{c}\end{equation}
As an intermediate step, we   prove the
following analogue of (\ref{DESIRED}):
\begin{lem} \label{R} Let $\Omega \subset \RR^n$ be a bounded piecewise smooth domain, with ergodic billiard map, and let $A_h$ be a zeroth order operator. For all $\epsilon > 0$, there exists a pseudodifferential operator
$R_h$ of the form (\ref{c}) such that
\begin{equation}
\limsup_{\lambda \to \infty} \frac1{N(\lambda)} \sum_{\lambda_j <
\lambda}  \Big| \langle (A_{h_j} -  \omega_{\Neu}(A) \; R_{h_j}) u_j^b,
\hat u_j^b \rangle \Big|^2 < \epsilon. \label{want-small}
\end{equation}
\end{lem}

The proof  uses  Cauchy-Schwarz  to bound the left side of \eqref{want-small} by
\begin{equation}
\limsup_{\lambda \to \infty} \frac1{N(\lambda)} \sum_{\lambda_j
< \lambda}  \llangle (A - R \omega_{\Neu}(A))^2 u_j^b,
u_j^b \rrangle . \label{BR2-intro}\end{equation}
 We then use the invariance properties (\ref{FONE}) and
   Corollary (\ref{INV})  to
replace $A$ by the time-average $A_N$ defined by
$$
(A_N)_h = \frac1{N} \sum_{k=1}^N ((F_h^k)^* A_h F_h^k.
$$
After replacing $A$ by $A_N$, it follows from the local Weyl law of
Lemma~\ref{LWL} that the limit in \eqref{BR2-intro} is given by
\begin{equation}
\omega_{\Neu}((A_N - \omega_{\Neu}(A) R )^2) = \int_{B^*Y} \Big(
\sigma(A_N) - \omega_{\Neu}(A) \sigma(R) \Big)^2 d \mu_{\Neu}.
\label{BR3-intro}\end{equation}
By applying Lemma~\ref{Egorov} iteratively we see that
the the symbol $\sigma(A_N)$ of $A_N$ is $\gamma(q)$ times a
discrete average of $\gamma^{-1} \sigma(A)$ over iterates of the billiard
map. Hence, this converges to $\gamma(q) \times \int \gamma^{-1} \sigma(A)  \, d\sigma / \int 1 \, d\sigma$, which is equal to $c \,
\gamma(q) \omega_{\Neu}(A)$. This is approximately equal to the
symbol of $\omega_{\Neu}(A) R$. Thus, the integral \eqref{BR3-intro}
becomes small as $N \to \infty$, and  we can make
\eqref{want-small} arbitarily small by choosing a sufficiently
good operator $R$.
   We   deduce (\ref{DESIRED}) by considering \eqref{want-small} with $A$ replaced by $\omega_{\Neu}(A) \Id$,
   and subtracting
 this from \eqref{want-small}.

In this sketch, we implicitly assumed that Egorov's theorem and
other pseudodifferential manipulations were valid on all of $Y$.
However, they are only valid away from the corners and hence we
have only obtained an equidistribution law away from the corners.
But  Lemma~\ref{LWL}, with $A_h = \Id$, shows that eigenfunctions do not concentrate at the corners and therefore $\omega_{\Neu}$ is the full limit measure.
 The proof is similar in the case of the other boundary
 conditions.

Let us now describe the organization of this article. Although the
details of the arguments differ somewhat (and sometimes
significantly) as the boundary condition $B$ varies,   we
present complete details first in the case of the Neumann boundary
condition, and then describe the necessary modifications for the
other boundary conditions.  Robin boundary condition comes second
because it is a deformation of the Neumann boundary condition and
the details are quite similar. We then take up the Dirichlet boundary
condition, which is singular relative to Neumann. The $\Psi^1$-Robin
boundary condition provide a bridge between Neumann and Dirichlet
and require the most serious modifications; we present the details
for these conditions last. We shall also assume that our domain is convex until the Section~\ref{nonconvex}, where we adapt the argument to nonconvex domains.
 To guide the proof, we also start in Section~\ref{ENDO} with an overview
 of how boundary ergodicity differs in principle
from the general structure of proof of quantum ergodicity given in
\cite{Z}. In the appendix we summarize properties of heat kernels
needed for the proof.

We wish to thank Ben Andrews, Nicolas Burq, Rafe
Mazzeo, and Tom ter Elst for helpful discussions; Leonid Bunimovich and Maciej Wojtkowski for
informing us about the state of the art on Euclidean domains with
ergodic billiards; Maciej Zworski for advising us that the generalization from convex to general domains should not be hard and encouraging us to treat this case; 
and Alan McIntosh, Andreas Axelsson, Monique Dauge and Michael Taylor for background on Lipschitz domains.


\section{\label{ENDO}Quantum ergodicity of  endomorphisms}

Quantum ergodicity is concerned with  quantizations of classically
ergodic Hamiltonian systems. The typical example is the wave group
$U_t = e^{i t \sqrt{\Delta}}$ of a compact Riemannian manifold
$(M, g)$ without boundary, whose geodesic flow $G^t: S^*M \to
S^*M$ is ergodic on the unit co-sphere bundle with respect to
Liouville measure $d\mu_L$. Ergodicity  on the quantum level is a
relation between the {\it time average}
$$\langle A \rangle := \text{w-}\lim_{T \to \infty} \frac{1}{2T}
\int_{-T}^T U_t A U_t^* dt$$ of an observable $A \in \Psi^0(M)$
and the {\it space average}, defined to be the constant operator
$$\omega(A) \; I, \;\;\;\; \mbox{where}\;\; \omega(A) = \frac{1}{\vol(S^*M)}
\int_{S^*M} \sigma_A d\mu_L.$$ Here, $\Psi^0(M)$ is the space of
zeroth order pseudodifferential operators over $M$ and $\sigma_A$
is the principal symbol of $A$.

The system is said to be {\it quantum ergodic} if
\begin{equation} \label{QE} \langle A \rangle = \; \omega(A)\; I\; + \; K,
\;\;\mbox{with}\;\;\; \frac{1}{N(\lambda)} || \Pi_{\lambda} K
\Pi_{\lambda}||^2_{HS} \to 0, \end{equation} where $||\cdot||_{HS}$
is the Hilbert-Schmidt norm and where $\Pi_{\lambda}$ is the
spectral projection onto the span of eigenfunctions of
$\sqrt{\Delta}$ of eigenvalue $\leq \lambda.$ Thus, time and space
averages agree up to an operator $K$ which is negligible in the
semi-classical limit. By expressing all operators in terms of the
eigenfunctions of $U_t$, one sees that this formulation is
equivalent to (\ref{DESIRED}).

The quantum dynamics is thus  the C* dynamical system $(\R,
\alpha_t, \Psi^0)$  defined by the one-parameter group of
automorphisms
\begin{equation} \label{AUTO} \alpha_t(A) = U_t A U_t^*
\end{equation} on the norm closure of $\Psi_0. $ It
 is a general fact that $(\R, \alpha_t, \Psi^0)$ is quantum
 ergodic as long as $\omega$ is an ergodic state \cite{Z}. After
 applying the Schwarz inequality for states and the Egorov
 theorem, the proof is reduced to the mean $L^2$ ergodic theorem
 on the symbol level.

 This general fact does not apply immediately to the wave group
 $U_t$ of the Laplacian $\Delta_B$ on a bounded domain $\Omega \subset M$
 in a Riemannian manifold $(M, g)$,
  with one of the above boundary conditions
 $B$. The problem is that conjugation by the wave group does not
 quite
 define an automorphism of any natural algebra of
 pseudodifferential operators. Nevertheless, by using suitable
 approximations, one can  prove that (\ref{DESIRED}) holds for
 pseudodifferential operators $A$ on the ambient manifold which
 are essentially supported in the interior of $\Omega$
 \cite{ZZw}.

 The situation we consider in this paper is complicated in one
 further way: As mentioned in the introduction, the relevant
 quantum dynamics (\ref{QDF}) defines an endomorphism rather than
 an automorphism of the relevant algebra of observables. Let us
first explain this on the classical level. For simplicity, we
restrict to smooth  plane domains and Neumann boundary conditions.
The general case is similar.

Let $\Omega \subset \R^n$ be a smooth domain and let
\begin{equation} \beta: B^* Y \to B^* Y \end{equation} denote the billiard map (see Section~\ref{corners}).
In the spectral theory of dynamical systems, one studies the
dynamics of $\beta$ through the associated `Koopman' operator
$$\mathcal{U}: L^2(B^*Y, d \sigma)
\to L^2(B^* Y, d\sigma),\;\; \mathcal{U} f (\zeta) =  f(\beta
(\zeta)).$$
 Here,
$d\sigma$ denotes the usual  $\beta$-invariant symplectic volume measure on
$B^* Y$. From the invariance it follows that $\mathcal{U}$ is a
unitary operator.  When $\beta$ is ergodic, the unique invariant
$L^2$-normalized eigenfunction is a constant $c$, and one has the
mean ergodic theorem (see e.g. \cite{Petersen}, chapter2)
\begin{equation} \label{UNITARY} \lim_{N \to \infty} ||\frac{1}{N}
\sum_{n = 1}^N \mathcal{U}^n(f) - \langle f, c \rangle || \to 0.
\end{equation}

Suppose now that we introduce a positive function  $\gamma \in
C^{\infty}(B^*Y)$, and define a new operator  $T$ by
$$T f (\zeta) = \frac{\gamma (\zeta)}{(\gamma(\beta(\zeta)))} f(\beta(\zeta)).$$ Then $T$ is not  unitary
on $L^2(B^*Y, d \sigma)$. This is the situation we find ourselves
in by virtue of the Egorov Lemma (\ref{Egorov}) in the case of
Neumann boundary conditions. We find that the orthogonal
projection $P$ onto the invariant $L^2$ functions for $T$ relative
to the invariant inner product  has a somewhat different form:

\begin{prop} Consider the operator $T$ on $L^2(B^*Y, d
\sigma).$ We claim:

\begin{itemize}
\item (i) The unique positive $T$-invariant $L^1$ function is given by  $\gamma.$
\item (ii)  $T$ is unitary relative to the   inner product $\langle \langle, \rangle \rangle$ on $B^* Y$
defined by the measure $d\nu = \gamma^{-2} d \sigma$.
\item (iii)  When $\beta$ is ergodic, the   orthogonal projection $P$ onto the invariant
$L^2$-eigenvectors has the form $$ P(f) = \frac{\langle \langle f,
\gamma \rangle \rangle}{ \langle \langle \gamma, \gamma \rangle
\rangle} \gamma = \frac1{\vol(B^*Y)} [\int_{B^* Y} f \gamma^{-1} d
\sigma] \; \gamma = c \omega_{\Neu}(f) \; \gamma
$$
where $c$ is as in \eqref{c}.
\end{itemize}

\end{prop}

\begin{proof} An
$L^1$ measure $\rho d \sigma$ is invariant under $T$ if and only if
$$\frac{\gamma (\zeta)}{\gamma (\beta(\zeta))} \rho (\beta(\zeta)) = \rho(\zeta) \iff
\gamma^{-1} (\zeta) \rho(\zeta) \;\; \mbox{is invariant under}
\;\; \beta.$$ The unique  positive solution (and the unique
solution in the ergodic case) is given by $ \rho d \sigma = \gamma d
\sigma. $ This proves (i).

To prove (ii), we  seek an invariant inner product $\langle
\langle, \rangle \rangle$ of the form
$$\langle \langle f, g \rangle \rangle = \int_B f(\zeta) \bar{g}(\zeta) d \nu. $$
We note that
$$\begin{aligned}\langle \langle T f, g \rangle \rangle & =   \int_B
\frac{\gamma ((\zeta)}{\gamma(\beta(\zeta))} f(\beta(\zeta)) g(\zeta) d\nu \\
& =  \int_B  f(\zeta) g(\beta^{-1}(\zeta)) \frac{\gamma (\beta^{-1}(\zeta))}{\gamma(\zeta)}  d\nu \circ \beta^{-1}  \\
& =  \llangle f, T^{-1} g \rrangle \iff   \frac{\gamma
(\beta^{-1}(\zeta)) }{\gamma(\zeta)}  d\nu \circ \beta^{-1} =
\frac{\gamma (\zeta)}{\gamma(\beta^{-1}(\zeta))}  d \nu \\
& \iff \gamma(\zeta)^{2} d \nu = \gamma(\beta(\zeta))^{2} d \nu \circ \beta \\
& \iff  d \nu = \gamma^{-2}  d \sigma.  \end{aligned}$$

The formula for $P$ follows from (i)-(ii).

\end{proof}

The formula for $P$ explains the need for the operator $R$ in
Lemma~\ref{R}. The  proof of quantum ergodicity will be based on a
reduction to the classical $L^2$-mean ergodic for the symbol using
Egorov's theorem.  In the case where $\beta$ is ergodic, the mean
ergodic theorem for $T$ on $L^2(B^* Y, \langle, \rangle)$ states
that:
$$\lim_{N \to \infty} ||\frac{1}{N} \sum_{n = 1}^N T^n(f) - P(f)|| \to 0. $$
In the automorphism case, $P(\sigma_A) = \omega(A)$ and we can
deduce the desired quantum ergodicity result directly from the
mean ergodic theorem.  In the endomorphism case, $P(\sigma_A) =
 c\omega_{\Neu}(A) \; \gamma$, so we need to construct an operator $R$
 which contributes the factor of $\gamma.$

 In the Dirichlet case, the roles of $T$ and $T^*$ get
 interchanged. By essentially the same proof, we then have:

 \begin{prop} Consider the operator $T^*$ on $L^2(B^*Y, d
\sigma).$ We claim:

\begin{itemize}

\item (i) The unique positive $T^*$-invariant density is given by $\gamma^{-1}.$

\item (ii)  $T^*$ is unitary relative to the   inner product $\langle \langle, \rangle \rangle$ on $B^*Y$
defined by the measure $d\nu = \gamma^{2} d \sigma$.

\item (iv)  When $\beta$ is ergodic, the orthogonal projection $P$ onto the invariant
$L^2$-eigenvectors has the form
$$
P(f) = \frac{\langle \langle f, \gamma^{-1} \rangle \rangle}{
\langle \langle \gamma^{-1}, \gamma^{-1} \rangle \rangle}
\gamma^{-1} = \frac1{\vol(B^*Y)} [ \!\!\! \int\limits_{B^* Y}
\!\!\! f \gamma d \sigma] \; \gamma^{-1} = c \omega_{\Dir}(f) \;
\gamma^{-1}.
$$
\end{itemize}

\end{prop}


\section{Piecewise smooth manifolds}\label{corners}
Let $\Omega$ be a bounded subdomain of $\RR^n$ with closure $\Omegabar$.

\begin{defn}\label{cpm} We say that $\Omega \subset \RR^n$ is a piecewise smooth manifold  if the boundary $Y = \partial \Omega$ is strongly Lipschitz, and can be written as a disjoint union
$$
Y = H_1 \cup \dots \cup H_m \cup \Sigma,
$$
where each $H_i$ is an open, relatively compact subset of a smooth  embedded hypersurface $S_i$, with $\Omega$ lying locally on one side of $H_i$, and where $\Sigma$ is a closed set of $(n-1)$-measure zero.

The sets $H_i$ are called
boundary hypersurfaces of $\Omega$. We call $\Sigma$ the singular set, and write $Y^o = Y \setminus \Sigma$ for the regular part of the boundary.
\end{defn}

\begin{rem} There appears to be no standard definition of `piecewise smooth manifold'.
 Our definition is rather broad; we have chosen it since it seems to include all known examples
 of ergodic Euclidean domains. We remark that it includes, for example, nonconvex polygons, and
 arbitrary convex polyhedra which are not included in the class of manifolds with corners as
 considered in \cite{CFS}, \cite{ZZw} and \cite{Burq}. However, the increase of generality is
  more apparent than real since the proofs in \cite{ZZw} and \cite{Burq} undoubtedly go through
   for piecewise smooth domains as defined here. Moreover, \cite{ZZw} and \cite{Burq} deal with
   Riemannian manifolds, not just Euclidean domains. 
\end{rem}

The operators $\Delta_B$ need to be properly defined on a piecewise smooth manifold. For all boundary conditions considered here, it is defined as the self-adjoint operator associated to a closed, semibounded quadratic form. The quadratic forms are
\begin{equation}
Q_{\Neu}(u) = \int_\Omega |\nabla u(z)|^2 \, dz \quad u \in H^1(\Omega),
\label{quadform-Neu}\end{equation}
for the Neumann boundary condition,
\begin{equation}
Q_{\operatorname{Rob}}(u) = \int_\Omega |\nabla u(z)|^2 \, dz   + \int_Y \kappa(y) |u(y)|^2 \, d\sigma(y) \quad u \in H^1(\Omega),
\label{quadform-Rob}\end{equation}
for the Robin boundary condition,
\begin{equation}
Q_{K}(u) = \int_\Omega |\nabla u(z)|^2 \, dz + \int_Y (Ku)(y) \overline{u(y)} \, d\sigma(y) \quad u \in H^1(\Omega),
\label{quadform}\end{equation}
for the $\Psi^1$-Robin boundary condition, and
\begin{equation}
Q_{\Dir}(u) = \int_\Omega |\nabla u(z)|^2 \, dz \quad u \in H^1_0(\Omega),
\label{quadform-D}\end{equation}
for the Dirichlet condition.
In \eqref{quadform}, $K$ is a first-order classical pseudodifferential operator whose kernel is supported in $Y^o \times Y^o$. Here closedness of \eqref{quadform} (hence self-adjointness of $\Delta_K$) requires that $K$ be self-adjoint, while ellipticity of the boundary problem requires that $K$ have non-negative principal symbol (the form is positive when $K$ is positive, and semi-bounded provided the principal symbol is non-negative).

\vskip 5pt

Suppose that $\Omega$ is a  piecewise smooth manifold.  We shall define a `compressed' cotangent space $\cTY = T^* Y^o \cup \Sigma$ which is the space of which our wavefront set, defined below, will be a subset. For each boundary hypersurface $H_i$ we define a compressed cotangent bundle $\tilde T^* H_i$ by identifying the fibre $T_y^* H_i$ above each $y \in \pa H_i$ to a point; we give this space the quotient topology. Then we define
the topology on $\cTY$ by gluing together the $\tilde T^* H_i$ along their boundaries. The space $\cTY$ essentially `ignores' the frequencies over $\Sigma$.

Let us now define certain classes of operators on $L^2(Y)$, depending on a
real parameter $h \in (0, h_0]$. First, we define $C^\infty(Y)$ to
be the restriction of $C^\infty(\RR^n)$ to $Y$, and define $C^\infty(Y^2)$ similarly. We define a {\it residual operator} to be one with a kernel $K_h(y,y')$ (all kernels are written with respect to Riemannian measure on $Y$) which are smooth in $(y,y')$ for sufficiently small $h$, and are $O(h^\infty)$ as an element of $C^\infty(Y^2)$ (in the sense that each seminorm of $K_h$ is $O(h^\infty)$).
We next define a {\it semiclassical pseudodifferential operator of order $m$} on $Y$ to be an operator which can be expressed as the sum of two terms, one a semiclassical pseudodifferential operator on $Y^o$, parametrized by $h \in [0, h_0]$ for some $h_0 > 0$,  whose kernel is properly supported in $Y^o \times Y^o$, and the second a residual operator. We shall work here only with the most classical types of operators, namely those which are either differential operators, or with symbols $a(y,\eta, h)$ which are smooth functions of $y, \eta$ and $h$, compactly supported in $\eta$.

 We similarly define a semiclassical Fourier integral operator (FIO) to be the sum of a semiclassical Fourier integral operator with kernel properly supported in $Y^o \times Y^o$, and a residual operator.
For a pseudodifferential operator or FIO $Q_h$, the (semiclassical) wavefront set support $\WF(Q_h)$ is well defined and is a closed subset of $T^* Y^o \times T^* Y^o$; for pseudodifferential operators, it is a subset of the diagonal.

We may define the operator wavefront set more generally.
Let $T_h$ be some operator on $Y$, and let $q, q' \in \cTY$. We say that $(q,q') \notin \WF (T)$ under the following conditions: if $q, q' \in \Sigma$, then if there are smooth functions $\phi$, $\phi'$ with $\phi(q) \neq 0$, $\phi(q') \neq 0$, and $\phi(y) T_h \phi(y')$ residual; if $q \in \Sigma$, $q' \in T^* Y^o$, then if there is $\phi$ as above and a pseudodifferential operator $Q_h$ which is elliptic at $q'$ with $\phi T_h Q_h$ residual; similarly if $q' \in \Sigma$, $q \in T^* Y^o$; and for $q, q' \in T^* Y^o$, if there are pseudodifferential  operators $Q, Q'$ elliptic at $q, q'$ respectively, with $Q_h T_h Q'_h$ residual. The last condition is equivalent to the usual condition when $Y$ is smooth. The operator wavefront set $\WF(T)$ is then a closed subset of $\cTY$.

\vskip 5pt

The billiard map $\beta : B^* Y^o \to \tilde T^* Y$ is defined on the open ball bundle $B^* Y^o$ as follows: given $(y, \eta) \in T^*Y$, with $|\eta| < 1$ we let $(y, \zeta) \in S^* \Omega$ be the unique inward-pointing unit covector at $y$ which projects to  $(y, \eta)$ under the map $T^*_{\pa \Omega} \Omegabar \to T^* Y$. Then we follow the geodesic (straight line) determined by $(y, \zeta)$ to the first place it intersects the boundary again; let $y' \in Y$ denote this first intersection.
If $y' \in \Sigma$ then we define $\beta(y, \eta) = y'$. Otherwise,  let $\eta'$ be the projection of $\zeta$ to $T^*_{y'}Y$. Then we define
$$
\beta(y, \eta) = (y', \eta').
$$
The map $\beta_{-} : B^* Y^o \to \tilde T^* Y$ is defined similarly, following the backward billiard trajectory
(that is, the straight line with initial condition $(y, 2 (\zeta \cdot \nu_{y}) \nu_{y} - \zeta )$).

The billiard map is a symplectic, hence measure preserving, map
with respect to the standard symplectic form on $T^*Y$. This follows from the fact that the
Euclidean distance function $d(y, y')$ is locally a generating function
for $\beta$; that is, the graph of $\beta$ in a neighbourhood of $(y_0, \eta_0, y'_0, \eta'_0)$ is given by
\begin{equation}
\{ (y, -\nabla_{y} d(y, y'), y',  \nabla_{y'} d(y,y') )
\label{local-d}\end{equation}
for $(y, y')$ in a neighbourhood of $(y_0, y_0')$.
We denote the graph of $\beta$ by $\Cbill$,
\begin{equation}
\Cbill = \operatorname{graph} \beta \equiv \{ (\beta(q), q) \mid q \in \mathcal{R}^1 \}.
\label{billiardCR}\end{equation}
This is a smooth Lagrangian submanifold of $B^* Y^o \times B^* Y^o$. For strictly convex $\Omega$ it is given globally by \eqref{local-d}, for $y, y' \in Y^o$,  but this is not true in general. This causes extra difficulties for nonconvex domains which are dealt with in Section~\ref{nonconvex}.

It is clear that when $\Sigma$ is nonempty,  $\beta$ does not map the open ball bundle to itself. We define
$\mathcal{R}^1 \subset B^* Y^o$ to be $\beta^{-1}(B^* Y^o)$ and, more generally, $\mathcal{R}^{k+1} = \beta^{-1}(\mathcal{R}^{k})$ for natural numbers $k$. Thus $\mathcal{R}^k$ consists of the points where $\beta^k$ is well defined and maps to $B^* Y^o$.
Similarly we define $\mathcal{R}^{-1} = \beta_-^{-1}(B^* Y^o)$ and $\mathcal{R}^{-k-1} = \beta_-^{-1}(\mathcal{R}^{-k})$.
Clearly $\mathcal{R}^1 \supset \mathcal{R}^2 \supset \dots$, $\mathcal{R}^{-1} \supset \mathcal{R}^{-2} \supset \dots$ and it is shown in \cite{CFS} that each $\mathcal{R}^k$ has full measure. Consequently, if we define $\mathcal{R}^\infty = \cap_k \mathcal{R}^k$, then $\beta$ is a  measure-preserving bijection on $\mathcal{R}^\infty$.


\section{Structure of the operators $E_h$ and $F_h$. }\label{Structure}
In this section, we shall decompose the operator $F_h$, defined by \eqref{Neumann-F}, as well as the related operator $E_h$, into microlocal pieces. Here $E_h$ is the operator with kernel
$$
E_h(y,y') = 2 G_0(y,y', h^{-1}) \Big|_{y, y' \in Y},
$$
where $G_0(y, y', \lambda)$ is the free outgoing resolvent kernel $(\Delta - (\lambda + i0)^2)^{-1}$ on $\RR^n$.  The first two pieces shall be standard types of operators, namely Fourier integral operators, which are well understood. The third piece is a `left over piece', but its operator wavefront set shall be quite localized.

The reader is reminded that in the section we assume that $\Omega$ is convex.

 For simplicity, we first consider the case when $Y$ is smooth. We denote the sphere bundle $\{ (y, \eta) \in T^*Y^o \mid |\eta| = 1 \}$ by $S^*Y^o$.

\begin{prop}\label{EandF} Assume that $\Omega$ is a smooth convex domain.  Let $U$ be any neighbourhood of $S^* Y^o \times S^* Y^o$. Then there is a decomposition of $E_h$
$$
E_h = E_{1,h} + E_{2,h} + E_{3,h},
$$
where $E_1$ is a Fourier Integral operator of order $-1$ associated with the canonical relation $\Cbill$ given by \eqref{billiardCR}, $E_2$ is a pseudodifferential operator of order $-1$ and $E_3$ has operator wavefront set contained in $U$. The principal symbol of $E_2$ in $B^*Y \setminus U_1$ (where $U_1$ is the projection of $U$ to the first factor) is $\-i/\gamma$, where $\gamma$ is defined by \eqref{a-defn}.

Similarly, there is a decomposition of $F_h$ as
$$
F_h = F_{1,h} + F_{2,h} + F_{3,h},
$$
where $F_1$ is a Fourier Integral operator of order zero associated with the canonical relation $\Cbill$, $F_2$ is a pseudodifferential operator of order $-1$ and $F_3$ has operator wavefront set contained in $U$.
\end{prop}

\begin{proof} We first tackle the simpler operator $E_h$.
The kernel of $E_h$ is given by
$$C \lambda^{n-2} (|y-y'|/h)^{-(n-2)/2}
\Ha_{n/2-1}(|y-y'|/h)
$$
The Hankel function $\Ha_{n/2-1}(t)$ is conormal at $t=0$ and as $t \to \infty$, $b(t) = e^{-it} \Ha_{n/2-1}(t)$ has an expansion in inverse powers of $t$, with leading term $\sim t^{-1/2}$. We introduce cutoff functions $1 = \phi_1(|y-y'|) + \phi_2(|y-y'|/h^{3/4}) + \phi_3(|y-y'|,h)$, where $\phi_1(t)$ is supported in $t \geq t_0$ for some $t_0 > 0$ to be chosen later, and $\phi_2(t)$ is equal to $1$ for $t \leq 1$ and equal to $0$ for $t \geq 2$.
(The power $3/4$ in $\phi_2$ could be replaced by any other power strictly between $1/2$ and $1$.)
Then, $\phi_1(|y-y'|) E_h(y,y')$ has a kernel of the form
$$
C h^{-(n-2)} e^{i|y-y'|/h} \phi_1(|y-y'|) \tilde b(|y-y'|/h),
$$
where $\tilde b(t)$ has an expansion in inverse powers of $t$ as $t \to \infty$, with leading term $\sim t^{-(n-1)/2}$. This is manifestly a semiclassical FIO of order $-1$, and since $\Omega$ is assumed convex, the phase function $|y-y'|$ parametrizes the billiard relation $\Cbill$.

We next show that the kernel $\phi_3 E_h$ has operator wavefront set supported in the set $U$, if $\delta$ is chosen sufficiently small.  It is sufficient to show that if $Q_h \in \Psi^0(Y)$ satisfies $\WF(Q_h) \circ U = \emptyset$, then $Q_h \circ (\phi_3 E_h)$ is residual, and similarly, if $\tilde Q_h \in \Psi^0(Y)$ satisfies $\WF(\tilde Q_h) \circ U = \emptyset$, then $(\phi_3 E_h) \circ \tilde Q_h$ is residual. The kernel of $Q_h \circ (\phi_3 E_h)$ is given by an oscillatory integral of the form
$$
h^{2-n} \int e^{i(y-y'') \cdot \eta/h} q(y'',\eta, h) e^{i|y''-y'|/h} \tilde b(|y''-y'|/h) \phi_3(|y''-y'|,h) \, dy'' \, d\eta
$$
where $q(y'', \eta)$ is supported away from $|\eta| = 1$. The phase is stationary when $\eta = d_{y''}|y'' - y|$. Since $\phi_3$ is supported in the region where $|y-y'| \leq \delta$, and $|d_{y''}(|y'' - y|)| \to 1$ as $|y-y'| \to 0$, this means that the phase is never stationary if $t_0$ is sufficiently small. Repeated integrations-by-parts show that the kernel is residual, since we gain an $h$ each time we differentiate the phase and lose at most $h^{3/4}$ when we differentiate $\phi_3$. The computation for $\tilde Q_h$ is similar.

Next we analyze the kernel $\phi_2 E_h$. We shall show that this kernel is pseudodifferential when microlocalized away from the set $|\eta| = 1$. To do this, we write the kernel of $\phi_2 E_h$ as  the distributional limit, as $\epsilon \to 0$, of the oscillatory integral
\begin{equation}
2(2\pi)^{-n} h^{2-n} \phi_2\Big( \frac{|y-y'|}{h^{3/4}} \Big) \int e^{i(y-y') \cdot \xi/h} \frac1{ |\xi|^2 - 1 - i\epsilon} \, d\xi .
\label{phi1-op}\end{equation}
We write $\xi = \xi^\parallel + \zeta \nu_{y'}$, where $\xi^\parallel \in T_{y'} Y$, $\zeta \in \RR$ and $\nu_{y'}$ is is the inward pointing unit normal at $y'$. The kernel can be written
\begin{equation*}
2(2\pi)^{-n} h^{2-n} \phi_2\Big( \frac{|y-y'|}{h^{3/4}} \Big) \int e^{i(y-y') \cdot \xi^\parallel /h}  \frac{e^{i(y-y') \cdot \nu_{y'} \zeta/h}}{\zeta^2 + |\xi^\parallel|^2  - 1 - i\epsilon} \, d\xi^\parallel \, d\zeta .
\end{equation*}
To localize away from $|\eta| = 1$, we introduce cutoffs $1 = \psi_1(\xi^\parallel) + \psi_2(\xi^\parallel) + \psi_3(\xi^\parallel)$, where $\psi_1(t)$ is supported  in $t \leq 1 - t_1$, $\psi_2(t)$ is supported in $1 - 2t_1 \leq t \leq 1 + 2t_1$ and $\psi_3(t)$ is supported in $t \geq 1 + t_1$. Inserting the cutoff $\psi_1$ means that $|\xi^\parallel|^2 - 1 < 0$, so we can perform the $\zeta$ integral using the formula
$$
\frac1{2\pi} \int e^{ik \zeta} \frac1{\zeta^2 - (a+i0)^2} \, d\zeta = \frac{i e^{ika}}{2a}, \quad k > 0, a > 0
$$
to get
\begin{equation}
i(2\pi)^{-n+1} h^{2-n} \phi_2\Big( \frac{|y-y'|}{h^{3/4}} \Big) \int e^{i(y-y') \cdot \xi^\parallel /h}  \frac{e^{i(y-y') \cdot \nu_{y'} \sqrt{1-|\xi^\parallel|^2 }/h}}{ \sqrt{1-|\xi^\parallel|^2 } }\psi_1(\xi^\parallel) \, d\xi^\parallel .
\label{phi1psi1}\end{equation}
As $h \to 0$, $(y-y') \cdot \nu_{y'}/h \to 0$ on the support of $\phi_2$, since $|y-y'| \leq 2 h^{3/4}$ on the support of $\phi_2$ but $(y-y') \cdot \nu_y = O(|y-y'|^2)$. Hence we may expand the exponential $e^{i(y-y') \cdot \nu_{y'} \sqrt{1-|\xi^\parallel|^2 }/h}$ in a Taylor series centred at zero:
\begin{multline}
\Big| e^{i(y-y') \cdot \nu_{y'} \sqrt{1-|\xi^\parallel|^2 }/h} - \sum_{j=0}^{N-1}
\frac1{j!} \Big( i(y-y') \cdot \nu_{y'} \sqrt{1-|\xi^\parallel|^2 }/h \Big)^j \Big| \\
\leq C_N \Big|(y-y') \cdot \nu_{y'} \sqrt{1-|\xi^\parallel|^2 }/h \Big|^N.
\label{phi1psi3}\end{multline}
Consider one of the terms
\begin{equation}
\frac{i}{j!} (2\pi)^{-n+1} h^{2-n}  \int e^{i(y-y') \cdot \xi^\parallel /h}
\Big( i(y-y') \cdot \nu_{y'} \sqrt{1-|\xi^\parallel|^2 }/h \Big)^j
\psi_1(\xi^\parallel) \, d\xi^\parallel .
\label{termj}\end{equation}
Since $(y-y') \cdot \nu_{y'} = O(|y-y'|^2)$, after $2j$ integrations by parts we eliminate the vanishing at the diagonal and gain a factor of $h^{2j}$. This is therefore a pseudodifferential operator of order $-j$. Multiplication by the $\phi_2$ factor only changes this by a residual kernel, since \eqref{termj} vanishes rapidly as $|y-y'|/h \to \infty$. Similarly, the error term is becoming more and more regular. Therefore \eqref{phi1psi3} is a pseudodifferential operator of order $-1$. Moreover, we see from the form of \eqref{phi1psi1} that the principal symbol of this operator is $i/\sqrt{1 - |\eta|^2} = i/\gamma$.

Similar reasoning applies to the cutoff $\psi_3$, using instead
$$
\frac1{2\pi} \int e^{ik \zeta} \frac1{\zeta^2 + a^2} \, d\zeta = \frac{e^{-ka}}{2a},
\quad k, \, a > 0.
$$
This gives the kernel
\begin{equation*}
(2\pi)^{-n+1} h^{2-n} \phi_2\Big( \frac{|y-y'|}{h^{3/4}} \Big) \int e^{i(y-y') \cdot \xi^\parallel /h}  \frac{e^{-(y-y') \cdot \nu_{y'} \sqrt{|\xi^\parallel|^2 - 1}/h}}{ \sqrt{|\xi^\parallel|^2 - 1} } \psi_1(\xi^\parallel) \, d\xi^\parallel .
\end{equation*}
We can similarly expand the exponential in a Taylor series to show that we get a pseudodifferential operator of order $-1$.

The operator \eqref{phi1-op} with cutoff $\psi_2$ inserted has operator  wavefront set arbitrarily close to $S^* Y^o \times S^* Y^o$, and hence within $U$ provided $t_1$ is sufficiently small. This is shown as for $\phi_2 E_h$, with the help of Theorem 7.7.1 of \cite{Ho} which gives $\epsilon$-independent estimates on all seminorms of the composition of the operator with $Q$ (on the left) or $\tilde Q$ (on the right).
Hence, if we define $E_1$ to be the operator $\phi_1 E$, $E_2$ to be the operator with cutoffs $\phi_2(\psi_1 + \psi_3)$ and $E_3$ to be the remainder, we have a decomposition which satisfies the conditions of the theorem.

\

To deal with the operator $F_h$, we argue similarly. Using the same $\phi_i$ cutoffs as before, the operator $\phi_1 F_h$ is an FIO of order $0$ and the operator $\phi_3 F_h$ has wavefront set contained in $U$. To deal with the remaining term, we write the kernel of $\phi_2 F_h$ as  the distributional limit, as $\epsilon \to 0$, of
\begin{equation}
2(2\pi)^{-n} h^{1-n} \int e^{i(y-y') \cdot \xi/h} \frac{-i\xi \cdot \nu_{y'}} { |\xi|^2 - 1 - i\epsilon} \phi_2\Big( \frac{|y-y'|}{h^{3/4}} \Big) \, d\xi .
\end{equation}
We decompose $\xi = \xi^\parallel + \zeta \nu_{y'}$ as before, and write the kernel
\begin{equation}
2(2\pi)^{-n} h^{1-n} \int e^{i(y-y') \cdot \xi^\parallel /h} e^{i(y-y') \cdot \nu_{y'} \zeta/h} \frac{-i\zeta}{\zeta^2 + |\xi^\parallel|^2  - 1 - i\epsilon} \phi_2\Big( \frac{|y-y'|}{h^{3/4}} \Big) \, d\xi .
\end{equation}
We use the cutoffs $1 = \psi_1(\xi^\parallel) + \psi_2(\xi^\parallel) + \psi_3(\xi^\parallel)$ as above. Inserting the cutoff $\psi_1$ means that $|\xi^\parallel|^2 - 1 > 0$, so we can perform the $\zeta$ (oscillatory) integral using
$$
\frac1{2\pi i} \int e^{ik \zeta} \frac{\zeta}{\zeta^2 - (a+i0)^2} \, d\zeta = \frac{e^{ika}}{2}, \quad k, \, a > 0
$$
to get
\begin{equation*}
i(2\pi)^{-n+1} h^{1-n} \int e^{i(y-y') \cdot \xi^\parallel /h} e^{i(y-y') \cdot \nu_y \sqrt{|\xi^\parallel|^2 - 1}/h} \phi_2\Big( \frac{|y-y'|}{h^{3/4}} \Big) \,
d\xi^\parallel .
\end{equation*}
Following the reasoning above, this appears to be a pseudodifferential operator with symbol $1$. Similarly, the term with $\psi_3$ appears to be a pseudo with symbol $1$, which would give us the identity operator modulo an operator of order $-1$. However, the identity term
is not present in the kernel $F_h$ since the kernel of the identity is supported at the diagonal and does not appear in the restriction of the kernel of $\pa_{\nu_{y'}} G_0(y, y')$ to the boundary. In fact, it is the `same' identity operator that turns up in the jump formula for the double layer potential in \eqref{jumpformula}. Thus, this piece of $F_h$ turns out to be of order $-1$, as for $E_h$.

Finally, as for $E_h$, the operator with cutoff $\psi_2$ may be shown to have wavefront set arbitrarily close to $\csd$, and hence within $U$ by a suitable choice of the cutoffs $\psi_i$. Hence, if we define $F_1$ to be the operator $\phi_1 F$, $F_2$ to be the operator with cutoffs $\phi_2 (\psi_1 + \psi_3)$ and $F_3$ to be the remainder, we have a decomposition which satisfies the conditions of the theorem.
\end{proof}

\begin{rem}
 A related analysis of $F_h$ is given in \cite{Z2}
(it is denoted $N(k + i \tau)$ there).
\end{rem}

A more complicated version of this Proposition is valid when $\partial \Omega$ has singularities. We first define
\begin{equation}
 \Xi = S^* Y^o \cup \Sigma \subset \tilde T^* Y,
\label{Xi}\end{equation}
and recall the notation $\mathcal{R}^k$ from Section~\ref{corners}. We denote the complement of a set $S$ by $S'$.

\begin{prop}\label{F-corners} Let $\Omega$ be a convex, piecewise smooth domain. Let $U$ be any neighbourhood of $\Xi \times (\mathcal{R}^1)' \cup (\mathcal{R}^{-1})' \times \Xi$. Then there is a decomposition of $F_h$
$$
F_h = F_{1,h} + F_{2,h} + F_{3,h},
$$
where $F_1$ is a Fourier Integral operator of order zero associated with the canonical relation $\Cbill$, $F_2$ is a pseudodifferential operator of order $-1$ and $F_3$ has operator wavefront set contained in $U$.
\end{prop}

\begin{proof}
We choose a function $\phi$ on $Y$ so that $\phi \equiv 1$ on a neighbourhood of $\Sigma$, so that
$$ \begin{aligned}
\{ (q, q') \mid \pi(q), \pi(q') \in \supp \phi \} &\subset U, \\
\{ (\beta(q), q) \mid \pi(q) \in \supp \phi \} &\subset U, \\
\{ (q, \beta^{-1}(q) ) \mid \pi(q) \in \supp \phi \} &\subset U.
\end{aligned}$$
We may write
$$
F = \phi F \phi + \phi F (1 - \phi)  + (1 - \phi) F \phi + (1 - \phi) F (1 - \phi).
$$
The term $(1 - \phi) F (1 - \phi)$ is supported away from the singular set in both variables, so may be treated by the argument above. We claim that the remaining terms have wavefront set contained in $U$. This is clear for $\phi F \phi$, and the argument for the other two terms is similar, so we concentrate just on $(1 - \phi)  F \phi$. We will show that the wavefront set is contained in
\begin{equation}
\{ (q, q') \mid \pi(q), \pi(q') \in \supp \phi \}  \cup
\{ (\beta(q), q) \mid \pi(q) \in \supp \phi \}
\label{wf-set}\end{equation}
Thus, let $(q,q')$ be a point not contained in \eqref{wf-set}. We want to show that the kernel $(1 - \phi)  F \phi$ is regular at $(q,q')$. If $\pi(q') \notin \supp \phi$, this is obvious, so assume that $\pi(q') \in \supp \phi$, $\pi(q') \notin \supp \phi$, and that $q = (y_0, \eta_0) \neq \beta(q')$. Then, there is a smooth function $a(y, \eta)$ with $\gamma(q) = 1$, a smooth function $\tilde \phi(y)$ with $\tilde \phi(\pi(q')) = 1$, and with the support of $a$ disjoint from
$$
\{ \beta(q') \mid \pi(q') \in \supp \tilde \phi \}.
$$
Let $A$ be a semiclassical pseudodifferential operator with symbol $a$.
Then the composition $A (1 - \phi) F \phi \tilde \phi$ is represented, modulo a residual term, by an integral of the form
\begin{equation}
\int e^{i(y-y'')\cdot \eta/h} e^{i|y''-y'|/h} a(y'', \eta) r(y'',y) (\phi \tilde \phi)(y') \, d\eta \, dy''.
\end{equation}
Here $r$ is a smooth function, since on the support of $a \tilde \phi$, $y''$ and $y'$ are separated.
The phase in this integral is never stationary on the support of the symbol, by construction. Hence this operator is residual, proving that $(q,q') \notin \WF( (1 - \phi) F \phi )$.
\end{proof}

In order to deal with products of operators involving $F_3$ or $F_3^*$ we need the following proposition.

\begin{prop}\label{WF}
Assume that the semiclassical operator $A_h$ is either a pseudodifferential operator or a Fourier integral operator with compact operator wavefront set. If $\WF(A) \circ \WF(F_3) = \emptyset$, then $AF_3$ is residual. Similarly, if $\WF(F_3^*) \circ \WF(A) = \emptyset$, then $F_3^* A$ is residual.
\end{prop}

\begin{proof}
We first observe that if $R_h$ is a residual operator, then $\WF(R F_3)$ is contained in $\{ (y, y) \mid y \in \Sigma \}$. Indeed, this property for $F$ follows readily from the form \eqref{Neumann-F}  of the operator, and it holds for $F_1$ and $F_2$ since these are FIOs. Therefore it holds also for $F_3 = F - F_1 - F_2$.

To prove that $A F_3$ is residual, notice that the sets
\begin{equation}\begin{aligned}
D_1 &= \{ q' \in \cTY \mid \ \exists (q, q') \in \WF(A) \}, \\
D_2 &= \{ q' \in \cTY \mid \ \exists (q', q'') \in \WF(F_3) \}.
\end{aligned}\end{equation}
are closed, and they are disjoint by hypothesis. Hence there exist disjoint open sets $O_1 \supset D_1$, $O_2 \supset D_2$. Choose $Q_h \in \Psi^0(Y)$ such that the symbol of $Q$ is $1$ in $O_1 \times (0, h_0)$ and zero in $O_2 \times (0, h_0)$ for some $h_0 > 0$, and so that the kernel of $Q$ is supported away from $\Sigma$ in both variables. Then $Q F_3$ is residual by hypothesis and and $A (\Id - Q) = A - AQ$ is residual by the symbol calculus for pseudodifferential operators, and supported away from $\Sigma$ in both variables. Writing
$$
A F_3 = A (Q F_3) + (A (\Id - Q)) F_3,
$$
and using the observation above we see that $AF_3$ is residual. The result for $F_3^* A$ follows by taking adjoints.
\end{proof}

We conclude this section with a crude operator bound on the operator $F_h$.

\begin{prop}\label{op-norm} Let $\Omega$ be a convex domain with
corners. Then  the $L^2$ operator norm of $F_h$ admits a bound
\begin{equation}
\| F_h \|_{L^2(Y) \to L^2(Y)} \leq C h^{-(n-1)}.
\label{op-bound}\end{equation}
\end{prop}

\begin{proof} The kernel of $F_h(z,z')$ is given by
$$
h^{-n} \nu_y \cdot (z - z') f( |z-z'|/h),
$$
where $f(t)$ is symbolic as $t \to 0$, with $f(t) = c t^{-n} + g(t)$, where $g(t)  = O(t^{-(n-1)})$ as $t \to 0$, and bounded as $t \to \infty$. Thus, the kernel is given by
$$
\nu_y \cdot \frac{z-z'}{|z-z'|^n} + h^{-(n-1)} \nu_y \cdot \frac{(z-z')}{h} g(|z-z'|/h).
$$
The first kernel is bounded on $L^2$ by the theory of singular integrals on Lipschitz submanifolds; see \cite{LMS}. The $L^2$ operator norm of the second kernel, whose singularity at the diagonal is $O(|z-z'|^{-n+2})$, may be crudely bounded by Schur's Lemma, giving a bound $C h^{-(n-1)}$.
\end{proof}

We remark that under any decomposition as in Proposition~\ref{F-corners},
 $F_1$ and $F_2$ are uniformly bounded on $L^2$, so the bound in \eqref{op-bound} is also valid for
 $F_3$.


\section{Local Weyl law}\label{Weyl}
In this section we shall prove Lemma~\ref{LWL} for the Neumann boundary condition.

Let us first prove \eqref{Weyl-limit} for $A_h = \Id$, which is the statement that
\begin{equation}
\lim_{\lambda \to \infty} \frac1{N(\lambda)}
\sum_{\lambda_j \leq \lambda} \langle u_j^b, u_j^b \rangle
\to \frac{2 \vol(Y)}{\vol(\Omega)}.
\label{psi-asympt}\end{equation}

We use the Karamata Tauberian theorem, which states (see eg \cite{T2}, p89)

\begin{prop}(Karamata Tauberian Theorem)\label{Kar} If $\mu$ is a positive measure on $[0, \infty)$, and $\alpha > 0$, then
$$
\int_0^\infty e^{-t\lambda} d\mu(\lambda) \sim a t^{-\alpha}, \quad t \to 0
$$
implies
$$
\int_0^x d\mu(\lambda) \sim \frac{a}{\Gamma(\alpha + 1)} x^\alpha, \quad x \to \infty.
$$
\end{prop}

To prove \eqref{psi-asympt} we consider
$$
e(t) = \sum_{j=1}^\infty e^{-t\lambda_j^2}  \langle u_j^b, u_j^b \rangle .
$$
This is equal to the trace of the operator $E(t)$ from \eqref{Ekernel} in the appendix. By \eqref{Id-trace}, we see that
$$e(t) \sim 2 (4\pi t)^{-n/2} \vol(Y) \text{ as } t \to 0,
$$
so by Proposition~\ref{Kar},
$$
\sum_{\lambda_j < \lambda} \langle u_j^b, u_j^b \rangle_Y \sim \frac{2 \vol(Y)}{(4\pi)^{n/2} \Gamma(n/2 + 1)} \lambda^n.
$$
On the other hand,
$$N(\lambda) \sim
\frac{\vol(\Omega)}{(4\pi)^{n/2} \Gamma(n/2 + 1)} \lambda^n,
$$
so \eqref{psi-asympt} follows.

Next we prove the lemma for zeroth order pseudos $A_h$  supported away from $\Sigma$
 in both variables. We proceed through a series of reductions.

(i) First, if $A_h$ is in $\Psi^{*,-1}(Y)$ (that is, $A_h = h \tilde A_h$ for $A_h \in \Psi^{*,0}(Y)$), then
\begin{equation}
\lim_{\lambda \to \infty} \frac1{N(\lambda)}
\sum_{\lambda_j \leq \lambda}  \langle A_{h_j}
u_j^b, u_j^b \rangle  \to 0.
\label{LWL(i)}\end{equation}
This follows easily from \eqref{psi-asympt} and the fact that the
operator norm of $\tilde A_h$ is uniformly bounded in $h$ (\cite{DSj}, Theorem 7.11).

(ii) Second, if the support of the symbol of $A$ at $h=0$ is contained in
$\{ |\xi| > 1 \}$, then we also have \eqref{LWL(i)}.

To see this, we use the fact that $F_{h_j} u_j^b = u_j^b$, to write
the left hand side of \eqref{LWL(i)} as
$$
\lim_{\lambda \to \infty} \frac1{N(\lambda)}
\sum_{\lambda_j \leq \lambda} \langle F_{h_j}^* A_{h_j}
F_{h_j} u_j^b, u_j^b \rangle.
$$
We use Proposition~\ref{F-corners} from the previous section to analyze
$F_h^* A_h F_h$. Due to the condition on the symbol of $A$, we can find an open set $U$ as in Proposition~\ref{F-corners} disjoint from $\WF(A)$.
 Let $F = F_1 + F_2 + F_3$ be a decomposition as in Proposition~\ref{F-corners} with respect to $U$. Then $\WF(A) \circ \WF(F_1) = \WF(A) \circ \WF(F_3) = \WF(F_1^*) \circ \WF(A) = \WF(F_3^*) \circ \WF(A) = \emptyset$. Thus by Proposition~\ref{WF}, the operators
$$
 A (F_1 + F_3) \text{ and } (F_1 + F_3)^* A
$$
are residual. Hence, using also Proposition~\ref{op-norm}, we have an $O(h^{\infty})$ estimate on the operator norm of $F_i^* A F_j$ unless $i = j = 2$, which certainly gives us the required estimate in these cases. In the remaining case, $i=j=2$, $F_2^* A F_2$ is a pseudodifferential operator of order $-2$, so the required estimate follows from (i).

(iii) Third, there exists an integer $k$ such that, if the symbol of $A_h$
at $h=0$ is sufficiently small in $C^k$ norm, then
\begin{equation}
\lim\sup_{\lambda \to \infty} \frac1{N(\lambda)}
\big| \sum_{\lambda_j \leq \lambda} \langle
A_{h_j}u_j^b, u_j^b \rangle \big| \leq \epsilon.
\label{LWL(iii)}\end{equation}

To prove this, note that by (i) one may assume that the symbol is
independent of $h$ and by (ii), that the support of the symbol is in $
\{ |\xi| \leq 2 \}$. Then, the $C^k$ norm of the symbol, for sufficiently
  large $k$, controls the operator norm of $A_h$, uniformly in $h$. Thus,
  if the $C^k$ norm is sufficiently small, the right hand side is
  smaller than $\epsilon$.

(iv) Since the symbols of differential operators are dense in the $C^k$
topology in the space of symbols supported in $ \{ |\xi| \leq 2 \}$, it is
enough, by (iii), to prove the local Weyl law for differential
operators.

(v) Next we note that for odd order monomial differential operators $A$, as above, the limit is zero. To see this, note that
\begin{equation*}
\frac{1}{N(\lambda)} \sum_{\lambda_j \leq \lambda}
\langle A_{h_j} u_j^b, u_j^b \rangle =  \\
\frac{1}{N(\lambda)} \sum_{\lambda_j \leq \lambda}
\frac{1}{2} \langle \big( A_{h_j} + A^*_{h_j} \big) u_j^b, u_j^b \rangle,
\end{equation*}
but $A + A^*$ is in $\Psi^{*,-1}(Y)$ if $A$ is an odd-order monomial, so by (i), the limit is zero. This clearly agrees with \eqref{Weyl-limit}, so it remains to treat even-order operators.

(vi) Finally, every even-order differential operator is the difference of two positive differential operators of the same order.
Hence we may restrict attention to positive even-order operators supported in a single coordinate patch.

\

Thus, let $A_h = h^{2k} P$, where $P$ is a positive differential operator of order $2k$. Consider the quantity
$$
e_A(t) = \sum_{j=1}^\infty e^{-t\lambda_j^2} \lambda_j^{-2k} \langle  P u_j^b, u_j^b \rangle .
$$
Let us consider $d_t^{k} e_A(t)$, which is given by
$$
\big( \frac{d}{dt} \big)^{k}e_A(t) = \sum_{j=1}^\infty e^{-t\lambda_j^2} \langle  P u_j^b, u_j^b \rangle  = \tr  P E(t).
$$
By \eqref{P-trace}, this is given by
\begin{equation}
\tr  P E(t) = \frac{4 (\frac{n}{2} + k-1)\dots (\frac{n}{2})  t^{-\frac{n}{2}-k}}{(4\pi)^{n/2} \vol(S^{n-1})}
\int\limits_{B^* Y} \sigma(P) d\mu_{\Neu} + O(t^{-\frac{n}{2}-k+\frac1{2}}).
\label{P-trace2}\end{equation}
Integrating $k$ times in $t$, we find
$$
e_A(t) =  \frac1{(4\pi t)^{n/2}} \frac{4}{\vol(S^{n-1})} \int\limits_{B^* Y} \sigma(P) d\mu_{\Neu} + O(t^{-n/2+1/2}), \quad t \to 0.
$$
Then \eqref{Weyl-limit} follows from Proposition~\ref{Kar} as before. This completes the proof of the Lemma.


\section{Egorov Theorem}\label{Egg}

In this section we prove a generalization of Lemma~\ref{Egorov}. First, we need to compute the principal symbol of the FIO $F_1$ from Proposition~\ref{F-corners}.

\begin{prop}\label{Eg} Let $F_1$ and $U$ be as in Proposition~\ref{F-corners}. Then the principal symbol of $F_1$ at $(\beta(q), q) \in \Cbill \setminus U$ is
$$
\sigma(F_1(\beta(q),q)) = \tau \Big( \frac{\gamma(q)}{\gamma(\beta(q))} \Big)^{1/2} |dq|^{1/2},
$$
where $\tau$ is an eighth root of unity and $|dq|$ represents the symplectic density on $B^* Y$.
\end{prop}

{\it Remark.} The Maslov bundle over $\Cbill$ is canonically trivial, since $\Cbill$ is the graph of a canonical transformation, so we shall ignore it.

\begin{proof}
We begin with the explicit expression \eqref{Neumann-F} for the free resolvent on $\RR^n$. At the billiard Lagrangian, $y$ and $y'$ are distinct, so we need the asymptotics of the Hankel function as it argument tends to infinity. This is
$$
\Ha_{n/2 - 1}(t) \sim e^{-i(n-1)\pi/4} e^{it} \sum_{j=0}^\infty a_j t^{-1/2 - j},
$$
where $a_0 = \sqrt{2/\pi}$. Moreover, the expansion can be differentiated term by term. Thus, the principal symbol of $F(\lambda)$ at the billiard Lagrangian (up to an eighth root of unity) is the same as that of
$$
(2\pi)^{-(n-1)/2} \lambda^{(n+1)/2} e^{i\lambda |y - y'|} |y - y'|^{-(n-1)/2} d_{\nu_{y'}} |y - y'| .
$$
We see from this expression that the
operator has semiclassical order zero, and its symbol is \cite{Ho4}, \cite{PUr}
$$
|y - y'|^{-(n-1)/2} d_{\nu_{y'}} |y - y'| \big| dy dy' \big|^{1/2}.
$$
It is more geometric to give the symbol in terms of the symplectic half-density $|dy d\eta|^{1/2} =|dy' d\eta'|^{1/2} = |dq|^{1/2} $ on $T^*Y$. To do this, we express $dy'$ in terms of $d\eta$, keeping $y$ fixed. Since $\eta_i = d_{y_i} |y - y'|$, we find that
$$
|d\eta| = \det \Big( \frac{\pa^2}{\pa s^i \pa t^j} \big| (y + s_i e_i) - (y' + t_i e'_i) \big| \Big)  |dy'|,
$$
where $e_i$ is an orthonormal basis for $T_y Y$, and $e'_i$ an orthonormal basis for $T_{y'}Y$. To make things a bit clearer we first compute this in the two dimensional case, $n=2$. We choose coordinates so that $y$ is at the origin, $y' = (0,r)$, $e_1 = (\cos \alpha, \sin \alpha)$ and $e_2 = (\cos \beta, \sin \beta)$. Thus, we are trying to compute
$$
\frac{\pa^2}{\pa_s \pa_t} \big| (s \cos \alpha - t \cos \beta, s \sin \alpha - r - t \sin \beta) \big|.
$$
This is equal to
$$
r^{-1} \cos \alpha \cos \beta = |y - y'|^{-1} \pa_{\nu_y} | y - y'| \pa_{\nu_{y'}} | y - y'|,
$$
so this gives the symbol
\begin{multline}
(d_{\nu_{y'}} |y - y'|) |y - y'|^{1/2} \Big( |y - y'|^{-1} \pa_{\nu_y} | y - y'| \pa_{\nu_{y'}} | y - y'| \Big)^{-1/2} \big| dy d\eta \big|^{1/2} \\
= \Big( \frac{d_{\nu_{y'}} |y - y'|}{d_{\nu_{y}} |y - y'|} \Big)^{1/2} \big| dy d\eta \big|^{1/2}  = \Big( \frac{\gamma(y', \eta')}{\gamma(\beta(y',\eta'))} \Big)^{1/2} \big| dy d\eta \big|^{1/2} ,
\label{F-symbol}\end{multline}
where $a$ is defined in \eqref{a-defn}.

In the higher dimensional case we compute $|d\eta|$ as follows:
let $X \subset \RR^n$ be the subspace of $\RR^n$
$$
X = (T_y Y) \cap (T_{y'} Y) \cap l_{y,y'}^{\perp},
$$
where $l_{y,y'}$ is the line joining $y$ and $y'$. We will initially assume that all intersections are transverse, so that $X$ has codimension three. Then $T_y Y \cap X$ has codimension two in $T_y Y$, and similarly for $T_{y'} Y$. We choose $e'_3 = e_3, \dots, e'_{n-1} = e_{n-1}$ to be an orthonormal basis of $X$.  We choose $e_1 \in T_y Y$ to be in the plane of $l_{y,y'}$ and $\nu_y$, and $e_2$ to be orthogonal to both $e_1$ and $l_{y,y'}$, and choose $e'_1$ and $e'_2$ similarly. Then, $e_1, \dots, e_{n-1}$ is an orthonormal basis of $T_yY$, and $e'_1, \dots, e'_{n-1}$ is an orthonormal basis of $T_{y'}Y$.

Let us choose coordinates in $\RR^n_z$ so that $X = \{ z_1 = z_2 = z_3 = 0 \}$, with $y= (0, \dots, 0)$, $y' = (0,0,r,0, \dots, 0)$. Without loss of generality  we may assume
$$\begin{aligned}
e_1 &= (\cos \alpha, 0 , \sin \alpha), \\
 e_2 &= (0,1,0), \\
 e'_1 &= (\cos \gamma \cos \beta, \cos \gamma \sin \beta, \sin \gamma), \\
e'_2 &= (-\sin \beta, \cos \beta, 0),
\end{aligned}$$
where we write only the first three coordinates since the rest are zeroes, and $e_j = e'_j$ is the $j+1$st unit vector in $\RR^n$. Thus we need to compute
\begin{equation}\begin{gathered}
\det  \frac{\pa^2}{\pa s^i \pa t^j} \big(  f(s_1, \dots, s_{n-1}, t_1, \dots, t_{n-1}) \big) , \quad \text{ where }\\
f = \big| (f_1, \dots, f_n) \big|, \quad \text{ and } \quad
\begin{aligned}
f_1 &= s_1 \cos \alpha - t_1 \cos \gamma \cos \beta + t_2 \sin \beta, \\
f_2 &= s_2 - t_1 \cos \gamma \sin \beta - t_2 \cos \beta, \\
f_3 &=  s_1 \sin \alpha - r - t_1 \sin \gamma,  \\
f_j &= s_{j-1} - t_{j-1}, \quad j \geq 4.
\end{aligned}
\end{gathered}\label{det-comp}\end{equation}
A simple but tedious computation shows that \eqref{det-comp} is equal to $r^{-n+1} \cos \alpha \cos \gamma$. If the transversality assumptions are not true, then this result follows by perturbing to a case where they are satisfied, and taking a limit. Thus, in general the symbol is given by \eqref{F-symbol}.
\end{proof}

{\it Remark.} The proof also shows that the principal symbol of $E_1$ at $(\beta(q), q) \in \Cbill \setminus U$ is equal to
\begin{equation}
\sigma(E_h)(\beta(q),q) = \frac{-i \tau}{\sqrt{\gamma(\beta(q)) \gamma(q)}} |dq|^{1/2}.
\label{E_1-symbol}\end{equation}

\begin{lem}\label{Egorov-N} Suppose $A_h$ is a pseudodifferential operator of order zero, with $\WF(A) \subset \mathcal{R}^{-N}$.
Then, one can express
$$
(F_h^*)^N A_h F_h^N = B_h + S_h,
$$
where $\| S_h \|_{L^2 \to L^2} \leq C h$, and $B_h$ is a pseudodifferential operator of order zero with symbol
\begin{equation}
 \sigma(B) = \begin{cases}
\gamma(q) \gamma(\beta^{N}(q))^{-1} a(\beta^{N}(q)) , \quad q \in \mathcal{R}^N \\
0 , \phantom{ \gamma(q)^{-1}\gamma(\beta(q)) b(\beta^{(N)}(q)) } \quad q \notin \mathcal{R}^N. \end{cases}
\label{Egorov-N-formula}\end{equation}
\end{lem}

\begin{proof}
Suppose that a neighbourhood $V$ of $(\Sigma \times \Sigma) \cup (\Sigma \times (\mathcal{R}^{-1})') \cup (\mathcal{R}' \times \Sigma) \cup \csd$ is given (we shall choose it later in the proof).
We decompose $F = F_1 + F_2 + F_3$ as in Proposition~\ref{F-corners}, with respect to $V$. Consider the product
\begin{equation}
(F_1^* + F_2^* + F_3^*)^N A (F_1 + F_2 + F_3)^N.
\label{product}\end{equation}

Consider terms in the expansion which contain at least one $F_3$ or $F_3^*$ (it will suffice to consider just $F_3$). These contain a string of the form
\begin{equation}
\dots A (F_1 \text{ or } F_2) (F_1 \text{ or } F_2) \dots (F_1 \text{ or } F_2) F_3 \dots
\label{string}\end{equation}
The substring $A (F_1 \text{ or } F_2) (F_1 \text{ or } F_2) \dots (F_1 \text{ or } F_2) = G$ is an FIO associated to the canonical transformation $\beta^k$, for some $k$ with $0 \leq k < N$. Thus, its operator wavefront set is contained in
\begin{equation}
\cup_{0 \leq k < N} \{ (q, \beta^{-k}(q)) \mid q \in \WF(A) \}.
\label{set}\end{equation}
By our assumption on the wavefront set of $A$, we have $\beta^{-k}(q)$  well defined for $q \in \WF(A)$. We now choose $V$ to be a neighbourhood of \begin{equation}
\Xi \times (\mathcal{R}^N)' \cup (\mathcal{R}^{-N})' \times \Xi
\label{set2}\end{equation}
such that no points of the form $(q, *)$  are in $V$, for $q \in \WF(A)$. This is possible since \eqref{set2} is compact, and the set $\WF(A) \times \cTY$ is a  closed set which is disjoint from \eqref{set2}. Then the FIO $G$ is such that $\WF(G) \circ \WF(F_3)$ is empty, so by Proposition~\ref{WF}, the composition \eqref{string} is a residual operator,  and hence is bounded on $L^2$ with an $O(h^\infty)$ bound. Since, by Proposition~\ref{op-norm}, we have a bound $C h^{-(n-1)}$ on all the remaining terms in the product \eqref{string}, it follows that each of these terms has an $O(h^\infty)$ bound on its operator norm, hence certainly are of the form $S_h$ above.

Thus, it remains to consider factors which do not contain an $F_3$ or $F_3^*$. Consider next all terms in the expansion of \eqref{product} which contain only $F_1$'s and $F_2$'s (including their adjoints), and at least one $F_2$. Since FIOs of order zero associated to a canonical transformation are bounded on $L^2$, and $F_2$ is order $-1$, all these terms are bounded on $L^2$ with an $O(h)$ bound, so these terms are of the form $S_h$ above.

We are thus reduced to analyzing the term $(F_1^*)^N B F_1^N$. We have shown that $F_1$ is an FIO whose canonical relation is $\Cbill = \operatorname{graph}(\beta)$ with symbol given by Proposition~\ref{Eg}, and hence, its adjoint $F_1^*$ is an FIO with canonical relation $\operatorname{graph}\beta^{-1}$, and with symbol
$$
\sigma(F_1^*)(q, \beta(q)) = \overline{ \sigma(F_1)(\beta(q), q) } =
\overline{\tau} \Big( \frac{\gamma(q)}{\gamma(\beta(q))} \Big)^{1/2} |dq|^{1/2}, \quad (\beta(q), q) \notin V.
$$
Thus, by the calculus of FIOs, the composition $(F^*_1) B F_1$ is a pseudodifferential operator with symbol
$\gamma(\beta(q))^{-1} \gamma(q) b(\beta(q))$.
By induction, we see that
the term $(F^*_1)^N B F_1^N$ is a pseudodifferential operator with symbol \eqref{Egorov-N-formula}.

 This completes the proof.

\end{proof}


\section{Proof of the main theorem --- Neumann boundary condition}\label{Proof}
In this section we prove \eqref{main-eqn} for the Neumann boundary condition.  It is sufficient to show that for any self-adjoint $A$,
\begin{equation}
\lim_{\lambda \to \infty} \frac1{N(\lambda)} \sum_{\lambda_j < \lambda} \big( \langle Au_j^b, u_j^b \rangle - \omega_{\Neu}(A) \big)^2 =0,
\end{equation}
(This gives us the limit \eqref{main-eqn} for a density one subsequence depending on $A$; a standard diagonal argument then shows that this dependence can be removed.)
As explained in
\S \ref{ENDO}, there are extra complications due to the fact that
\eqref{QDF} is not an automorphism. The main novelty in the
proof beyond consists in how to get around them.

The proof consists of several steps. First we show that
\begin{equation}
\lim_{\lambda \to \infty} \frac1{N(\lambda)} \sum_{\lambda_j < \lambda} \Big( \langle Au_j^b, \hat u_j^b \rangle - \frac{\omega_{\Neu}(A)}{\omega_{\Neu}(\Id)} \| u_j^b \|_{L^2(Y)} \Big)^2 = 0.
\label{step1-alt}\end{equation}
(Recall that $\hat u_j^b$ is the normalized boundary trace: $\hat u_j^b = u_j^b / \| u_j^b \|_{L^2(Y)}$.)
It is sufficient to show that for any $\epsilon > 0$, we have
\begin{equation}
\limsup_{\lambda \to \infty} \frac1{N(\lambda)} \sum_{\lambda_j < \lambda}  \Big( \langle Au_j^b, \hat u_j^b \rangle - \frac{\omega_{\Neu}(A)}{\omega_{\Neu}(\Id)} \| u_j^b \|_{L^2(Y)} \Big)^2 < C\epsilon,
\label{step1}\end{equation}
where $C$ depends only on $\Omega$ and $A$.
Notice that in \eqref{step1}, exactly {\it one} of the $u_j^b$'s is $L^2$-normalized. This $L^2$-normalization is a nuisance, but is not surprising, because in our proof of \eqref{step1} we make no reference to interior values of the eigenfunction, and the norm of $u_j^b$ in $L^2$ is dictated by the size of the eigenfunction in the interior.

In the second step, we show that
\begin{equation}
\| u_j^b \|^2_{L^2(Y)} \to \omega_{\Neu}(\Id)
\label{step2}\end{equation}
along a density one subsequence. Then \eqref{step1} and \eqref{step2} together imply the result of Theorem~\ref{main}.

\

{\it Step 1: Proof of \eqref{step1}.} We first choose a self-adjoint pseudodifferential operator $R$ whose principal symbol $r$ is a good approximation to $c \gamma$, in the sense that
\begin{equation}
\big\| \gamma^{-1/2} \big( r - c \gamma \big) \big\|_{L^2(B^*Y)}^2 < \epsilon, \quad c = \frac{\vol(S^{n-1}) \vol(\Omega)}{4 \vol(B^*Y)}.
\label{r-bound}\end{equation}
We choose a real symbol $a''$ such that $a''$ vanishes at the boundary of $B^*Y$, and we have
$$
\omega_{\Neu}(a'') = \omega_{\Neu}(a), \quad \omega_{\Neu}((a - a'')^2) < \epsilon.
$$
(Here we abuse notation by writing $\omega_{\Neu}(a)$ for $\omega_{\Neu}(A)$, etc, but this should not cause any confusion.)
We then define
\begin{equation}
{a''}^{(N)} = \frac1{N} \sum_{k=1}^N \gamma(q) \frac{ a''(\beta^k(q)) }{\gamma(\beta^k(q))},
\label{a''N}\end{equation}
an averaged version of the symbol in \eqref{Egorov-N-formula}.
By the mean ergodic theorem, we have
$$
\gamma^{-1}(q) {a''}^{(N)}(q) \to \frac1{\vol(B^* Y)} \int_{B^* Y} \gamma^{-1} {a''}^{(N)} d\sigma =
c \omega_{\Neu}(a'')  = c \omega_{\Neu}(A)\ \text{ in } L^2.
$$
(Note that we cannot make this claim for $a$ since $a \gamma^{-1}$ may not be $L^2$, due to the simple vanishing of $\gamma^2$ at the boundary of $B^*Y$; this is the sole reason for introducing $a''$. ) Here $c$ is as in \eqref{r-bound}.
We choose an integer $N$ such that
\begin{equation}
\big\| \gamma^{1/2} \big( \gamma^{-1}{a''}^{(N)} - c \omega_{\Neu}(A) \big) \big\|_{L^2(B^*Y)}^2 < \epsilon.
\end{equation}
We then choose a self-adjoint pseudodifferential operator $A'$ whose symbol $a'$ is close to $a$, and with additional properties required by Lemma~\ref{Egorov-N}. Namely, we require that $\WF(A')$ is contained in $\mathcal{R}^N \cap \mathcal{R}^{-N}$, and disjoint from $\Delta_{S^*Y}$, and is such that
\begin{equation}
\big( \omega_{\Neu}(A) - \omega_{\Neu}(A') \big)^2 < \epsilon, \quad
\omega_{\Neu}((a'' - a')^2) < \epsilon.
\label{A'}\end{equation}
Define ${a'}^{(N)}$ analogously to ${a''}^{(N)}$ in \eqref{a''N}. Then writing
$$
\gamma^{-1/2} \big( {a'}^{(N)} - c \gamma \omega_{\Neu}(A) \big) =
\gamma^{-1/2} \big( {a'}^{(N)} - {a''}^{(N)} \big) +
\gamma^{-1/2} \big( {a''}^{(N)} - c \gamma \omega_{\Neu}(A) \big),
$$
and applying the elementary inequality for real numbers
\begin{equation}
\big( a_1 + \dots + a_n \big)^2 \leq n \big( a_1^2 + \dots + a_n^2 \big)
\label{ineq}\end{equation}
with $n=2$, we find
\begin{equation}
\big\| \gamma^{-1/2} \big( {a'}^{(N)} - c \gamma \omega_{\Neu}(A) \big) \|_2^2 < 4\epsilon.
\label{a'N-bound}\end{equation}

\vskip 5pt

We then investigate the quantity
\begin{equation}
\lim \sup_{\lambda \to \infty} \frac1{N(\lambda)} \sum_{\lambda_j < \lambda}  \Big| \langle (A -  \omega_{\Neu}(A)R) u^b_j, \hat u^b_j \rangle \Big|^2 .
\label{BR}\end{equation}
To work with this, we want to replace $A$ with $A'$ since $A'$ has been constructed so that Lemma~\ref{Egorov-N} can be applied. Using \eqref{ineq} with $n=3$, we bound \eqref{BR} by
\begin{equation}\begin{gathered}
3 \lim \sup_{\lambda \to \infty} \frac1{N(\lambda)} \Bigg\{ \sum_{\lambda_j < \lambda}  \Big| \langle (A' -  \omega_{\Neu}(A')R) u^b_j, \hat u^b_j \rangle \Big|^2 \\
+ \sum_{\lambda_j < \lambda}  \Big| \langle (A -  A') u^b_j, \hat u^b_j \rangle \Big|^2
+ \sum_{\lambda_j < \lambda}  \Big| \langle (\omega_{\Neu}(A) -  \omega_{\Neu}(A'))R  u^b_j, \hat u^b_j \rangle \Big|^2 \Bigg\}.
\label{3eqn}\end{gathered}\end{equation}
In the first term of \eqref{3eqn}, we may replace $A'_h$ by
$$
{A_h'}^{(N)} = \frac1{N} \sum_{k=1}^N (F_{h}^*)^k A_h F_h^k,
$$
since the $u_j^b$ are invariant under $F_{h_j}$. Also using the inequality
$$
 \Big| \langle T u^b_j, \hat u^b_j \rangle \Big|^2 \leq \langle T^2 u_j^b, u_j^b \rangle
$$
for self-adjoint operators $T$, we can bound \eqref{3eqn} by
\begin{equation}\begin{gathered}
3 \lim \sup_{\lambda \to \infty} \frac1{N(\lambda)} \Bigg\{ \sum_{\lambda_j < \lambda}  \langle ({A'}^{(N)} -  \omega_{\Neu}(A')R)^2 u^b_j, u^b_j \rangle \\
+ \sum_{\lambda_j < \lambda}  \langle (A -  A')^2 u^b_j,  u^b_j \rangle
+ (\omega_{\Neu}(A) -  \omega_{\Neu}(A'))^2 \sum_{\lambda_j < \lambda}  \langle R^2  u^b_j, u^b_j \rangle  \Bigg\}.
\label{33eqn}\end{gathered}\end{equation}

We now use the Local Weyl Law, that is, Lemma~\ref{LWL}, to estimate the  three lim sups in \eqref{33eqn}. The second is bounded by $\epsilon$ by \eqref{A'}. The third is given by
$$
\frac{4 (\omega_{\Neu}(A) -  \omega_{\Neu}(A'))^2}{\vol(S^{n-1}) \vol(\Omega) } \int_{B^*Y} \frac{r(q)^2}{\gamma(q)} dq.
$$
Using \eqref{A'} and \eqref{ineq} with $n=2$ we can bound this by
$$
C \epsilon \int_{B^*Y} \frac{(r(q) - c \gamma(q))^2}{\gamma(q)} + \frac{( c \gamma(q))^2}{\gamma(q)}dq,
$$
where $C$ depends only on $\Omega$. Using \eqref{r-bound} we see that this is bounded by $C\epsilon^2 + C' \epsilon$.

To estimate the first term, we use Lemma~\ref{Egorov-N} and express
${A'}^{(N)} = \tilde A^{(N)} + S$, where $S$ has an $O(h)$ operator norm bound and $\tilde A^{(N)}$ is a pseudodifferential operator with symbol given by ${a'}^{(N)}$.
By \eqref{psi-asympt} we may neglect the operator $S$, and replace ${A'}^{(N)}$ by $\tilde A^{(N)}$.
Now, using the local Weyl Law for the operator $(\tilde A^{(N)} - \omega_{\Neu}(A') R)^2$, the first term is given by
\begin{equation}\begin{gathered}
\omega_{\Neu}(\tilde A^{(N)} - \omega_{\Neu}(A') R)^2 =
C\int_{B^*Y} \frac{ \big( {a'}^{(N)} - \omega_{\Neu}(A')r \big)^2 }{\gamma} \, dq
\\
\leq 2C \Bigg\{ \int_{B^*Y} \frac{ \big( {a'}^{(N)} - c \gamma \omega_{\Neu}(A') \big)^2 }{\gamma} \, dq
+ \omega_{\Neu}(A')^2 \int_{B^*Y} \frac{ \big( r - c \gamma \big)^2 }{\gamma} \, dq  \Bigg\}.
\end{gathered}\end{equation}
The first term is the key term which is controlled by \eqref{a'N-bound} thanks to classical ergodicity of the billiard flow.  The second term is controlled by \eqref{r-bound}.
This proves that \eqref{BR} is bounded by a constant times $\epsilon$.

\

Now \eqref{step1} follows by replacing $A$ in \eqref{BR} with $\omega_{\Neu}(A) \cdot \Id$, and then subtracting the two expressions:
\begin{equation}\begin{gathered}
\lim \sup_{\lambda \to \infty} \frac1{N(\lambda)} \sum_{\lambda_j < \lambda}  \Big| \langle (A -  \frac{\omega_{\Neu}(A)}{\omega_{\Neu}(\Id)}) u^b_j, \hat u^b_j \rangle \Big|^2  \\
\leq 2 \lim \sup_{\lambda \to \infty} \frac1{N(\lambda)} \sum_{\lambda_j < \lambda}  \Bigg\{ \Big| \langle (A -  \omega_{\Neu}(A)R) u^b_j, \hat u^b_j \rangle \Big|^2   \\ + \Big| \langle (\frac{\omega_{\Neu}(A)}{\omega_{\Neu}(\Id)} - \omega_{\Neu}(A)R ) u^b_j, \hat u^b_j \rangle \Big|^2  \Bigg\}
\end{gathered}\end{equation}
which we have shown is bounded by $C\epsilon$, where $C$ depends only on $\Omega$ and $A$. This completes the proof of \eqref{step1}.

\

{\it Step 2: Proof of \eqref{step2}. } To prove this, we need to relate the boundary values of the $u_j$ to their interior values, and use the $L^2$-normalization of the $u_j$ on $\Omega$. We do this via Green's formula:
\begin{equation}
u_j(z) = \int_Y G_{h_j}(z,y) u_j^b(y) \, d\sigma(y).
\label{Green}\end{equation}
Here, $G_h(z,y)$ is the normal derivative, in $y$, of the free Green function for $(-\Delta - (h^{-1} + i0)^2 )^{-1}$ on $\RR^n$. (Thus, $F_h$ is simply the restriction of this kernel to the boundary in $z$.) Let $\psi$ be a smooth nonnegative function compactly supported in the interior of $\Omega$. Then
\begin{equation}
\langle \psi^2 u_j,  u_j \rangle_{\Omega} =
 \langle G_{h_j}^* \psi^2 G_{h_j} u_j^b, u_j^b \rangle_Y .
\label{interior}\end{equation}

We need to analyze the operator $G_{h}^* \psi^2 G_{h}$.

\begin{lem}\label{G} Let $\phi$ be a smooth function on $Y = \cup_i H_i$ which is identically one near $\Sigma$. Then

(i) $G_{h}^* \psi^2 G_{h} (1 - \phi)$ is a semiclassical pseudodifferential operator of order zero on $Y$;  and

(ii) For each $1 \leq i \leq m$, $G_{h}^* \psi^2 G_{h} $ is bounded on $L^2(H_i)$, uniformly in $h$.

(iii) The limit of $\omega_{\Neu}(G^* \psi^2 G (1 - \phi))$, as the support of $\phi$ shrinks to $\Sigma$, is
\begin{equation}
\frac{1}{\vol(\Omega)} \int_\Omega \psi^2(z) \, dz.
\label{S-neu}\end{equation}
\end{lem}

\begin{proof}

(i) The kernel of the composition $(\psi G_h)^* (\psi G_h)(1 - \phi)$, as a function of $(y, y') \in Y \times Y$, is given by an integral of the form
\begin{equation}
h^{-(n-1)} (2\pi)^{-(n-1)} \int_\Omega e^{i |y-z|/h} e^{-i |y'-z|/h} a(y, z, h) \overline{a}(y', z, h) \psi^2(z) (1 - \phi)(y') \, dz.
\label{comp-int}\end{equation}
Because $y, y'$ are restricted to $Y$, and $z$ is restricted to the support of $\psi$ which is in the interior of $\Omega$, we do not see the singularity in the functions $| y - z|$, $|y' - z|$, or in $a$. The derivative of the phase function $\Phi(y, y', z) = | y - z| - | y' - z|$ in $z$ is nonvanishing unless $y = y'$. Due to the cutoff function $1 - \phi(y')$, which is supported in the regular part $Y^o$ of the boundary, the kernel is smooth and order $h^\infty$ unless $y = y'$ and both lie in a compact subset of $Y^o$.

To analyze further we write $\Phi(y, y', z) = \sum_i (y - y')_i \eta_i$, where $\eta_i(y, y', z)$ is a smooth function, which we may do since $\Phi$ vanishes whenever $y = y'$. Moreover, we have $\eta_i(y, y, z) = d_{y_i}\Phi(y, y', z)_{ | y' = y} = \Pi_y(y-z/|y-z|)$, where $\Pi_y$ is the projection onto the tangent space at $y$. Letting $s = |y - z|$, we find that
$$
s^{n-1} ds \, d\eta = \gamma(y, \eta) dz, \text{ when } y = y'.
$$
Hence we may change variable from $z$ to $(s,\eta)$ in a neighbourhood of the set $y = y'$, $y, y' \in \operatorname{supp} 1 - \phi$, $z \in \operatorname{supp} \psi$. There, the kernel takes the form
\begin{equation*}\begin{gathered}
h^{1-n} (2\pi)^{-(n-1)} \int e^{i(y-y') \cdot \eta} b(y, y', \eta, h) \, d\eta, \quad \text{where }
b \text{ is } C^\infty, \\
b(y, y, \eta, 0) = \int\limits_0^\infty a(y, z ,0), \overline{a}(y, z, 0)  \psi^2(z) (1 - \phi(y')) \, \frac{s^{n-1} ds}{ \gamma(y, \eta) }, \quad z = z(y, s, \eta),
\end{gathered}\end{equation*}
which is the kernel of a semiclassical pseudodifferential operator $B_h$ of order zero. (In the equation above, $z = y + s \theta$, where $\theta \in S^{n-1}$, $\theta \cdot \nu_y = \eta$, $\theta$ is inward pointing at $y$.) Moreover, we see that the symbol $b$ is nonnegative, and has the property that
\begin{equation}
\omega_{\Neu}(B_h) \ \text{ increases as either $\psi$ or $1 - \phi$ increases.}
\label{inc}\end{equation}.

(ii) This follows from (i).  Let $H_1, \dots, H_m$ be an enumeration of the  boundary hypersurfaces of $\Omega$. Each $H_i$ may be embedded in an open smooth submanifold $\tilde H_i$. Let $\chi$ be a function on $\tilde H_i$ with compact support, which is identically $1$ on $H_i$. Then, by (i), $G_{h}^* \psi^2 G_{h} \chi$ is a pseudodifferential operator of order zero, so is uniformly bounded. This implies a uniform bound $M$ on $\psi G_{h}$ acting from $L^2(H_i)$ to $L^2(\Omega)$. Then $m^2 M$ is a uniform bound for $G_{h}^* \psi^2 G_{h} \chi$ on $L^2(Y)$.

(iii) To prove this, we analyze the integral \eqref{comp-int} further. We are only interested in the principal symbol of $G^* \psi^2 G (1 - \phi)$, so we may substitute the value of $|a|$ at $h=0$, which is
$$
|a(y, z, 0)| = \frac1{2}  |y - z|^{-(n-1)/2} \gamma(y, z) \ \text{ since } \gamma(y, z) = \frac{y-z}{|y-z|} \cdot \nu_y.
$$
Thus, the principal symbol is given by
$$
b(y, y, \eta, 0) = \frac1{4} \int\limits_0^\infty s^{-(n-1)} \gamma(y, z)  \psi^2(z) (1 - \phi(y)) \, s^{n-1} ds,  \ z = z(y, s, \eta).
$$
Thus, $\omega_{\Neu}(G^* \psi^2 G (1 - \phi))$ is given by
$$
\frac{1}{\vol(S^{n-1}) \vol(\Omega)} \int\limits_{B^*Y} \int\limits_0^\infty  \psi^2(z) (1 - \phi(y)) \gamma(y, z) \, ds \frac{dy d\eta}{\gamma(y, \eta)} .
$$
The factors of $\gamma$ cancel.
Now we change variables of integration to $z$ and $\theta \in S^{n-1}$, where $z = y + s \theta$ and $\eta = \theta \cdot \nu_y$.  Then, $d\eta = \gamma d\theta$, and for a fixed $\theta$, $dz = \gamma dy \, ds$, so $dy \, ds \, d\eta = dz \, d\theta$ and the integral becomes
$$
\frac{1}{\vol(S^{n-1}) \vol(\Omega)} \int_{S^{n-1}} d\theta \int_{\Omega}  \psi^2(z) (1 - \phi(y(z, \theta))) \, dz.
$$
Clearly as the support of $\phi$ shrinks to the singular set, this integral converges to \eqref{S-neu}.
\end{proof}

We shall use the relation \eqref{interior}, Lemma~\ref{G} and the following result to deduce \eqref{step2}.

\begin{lem}\label{lem-seq}  Let $(a_n)$ be a seqence of complex numbers. Suppose that for every $\delta > 0$, there is a set of integers $S_\delta$ of density at least $1 - \delta$ such that, along it, the oscillation of the sequence $(a_n)_{n \in S_{\delta}}$ is at most $\delta$. Then there is a set $T$ of density $1$ such that the sequence $(a_n)_{n \in T}$ converges.
\end{lem}

{\it Remark. } The oscillation of a sequence $(b_n)$ is defined to be
$$
\lim_{N \to \infty} \ \sup_{m,n > N} | b_m - b_n |.
$$

\begin{proof}
By hypothesis, for each $n$ there is a set of integers $S_n$ of density at least $1 - 2^{-n}$ such that the oscillation of the corresponding subsequence is at most $2^{-n}$. By replacing $S_n$ by $S_n \cup \{ 1, 2, \dots, K \}$ for suitable $K$, we can ensure that
$$
\frac{ \# S \cap \{ 1, 2,  \dots M \} }{M} > 1 - 2^{-n+1} \ \text{ for all } M.
$$
Consequently, $T_n = S_n \cap S_{n+1} \cap \dots $ has density at least $1 - 2^{-n+2}$, and the sequence $(a_k)_{k \in T_n}$ has zero oscillation, that is, is a Cauchy sequence, hence converges to some number $A_n$. Since $T_n \cap T_m$ has positive density, for $m > n \geq 3$, we conclude that the $A_n$ are all equal to some fixed $A$.

By replacing $T_n$ by $T_n \setminus \{ 1, 2, \dots K \}$ for suitable $K$, we can ensure that $| a_k - A |  \leq 2^{-n}$ for all $k \in T_n$. Assuming this condition, then, the set $T = \cup_n T_n$ is density one, and $(a_n)_{n \in T}$ converges. To prove this, let $n$ be given. Then choose $L_i$ so that
$$
\sup | a_k - A | \leq 2^{-n} \ \text{ for } k \in T_i, \ k > L_i , \quad 1 \leq i \leq n.
$$
Let $L = \max_{i \leq n} L_i$; then
$$
\sup | a_k - A | \leq 2^{-n} \ \text{ for } k \in V, \ k > L , \quad 1 \leq i \leq n.
$$
This proves that the sequence $(a_n)$ converges to $A$ along $T$.
\end{proof}

Thus, by the Lemma, it suffices to show for each $\delta > 0$ that there is a sequence $S_\delta$ of density at least $1 - C\delta$ such that the oscillation in the sequence $\| u_j^b \|_{L^2(Y)}$ along $S_\delta$ is at most $C\delta$, for some constant $C$ depending only on $\Omega$.

Returning to \eqref{interior}, we write $G^* \psi^2 G = S$ and consider the equality
\begin{equation}
\langle \psi^2 u_j, u_j \rangle_{\Omega} = \| u_j^b \|_{L^2(Y)}^2 \langle S(1 - \phi) \hat u^b_j, \hat u^b_j  \rangle_Y
+ \langle S \phi u^b_j, u^b_j  \rangle_Y.
\label{uj-norm}\end{equation}
By Lemma~\ref{G}, the operator $S(1 - \phi)$ is pseudodifferential of order zero. Thus, by \eqref{step1}, we see that $\langle S(1 - \phi) \hat u^b_j, \hat u^b_j  \rangle_Y$ has a limit, namely $\vol(\Omega)(\vol(Y))^{-1}\omega_{\Neu}(S(1 - \phi)) > 0$, along a density one subsequence. Now we choose $\psi$ so that $0 \leq \psi \leq 1$, and $\psi = 1$ on a set of measure at least $(1 - \delta^2) \vol(\Omega)$. Then, by the local Weyl law in the interior of $\Omega$,
$$
\liminf_{\lambda \to \infty} \frac1{N(\lambda)} \sum_{\lambda_j < \lambda}   \langle \psi^2 u_j, u_j \rangle_{\Omega} \geq 1 - \delta^2.
$$
On the other hand, each individual term $\langle \psi^2 u_j, u_j \rangle_{\Omega}$ lies between $0$ and $1$, so we conclude that there is a subsequence $S'_\delta$ of density at least $1 - \delta$ on which the oscillation of the sequence $\langle \psi^2 u_j, u_j \rangle_{\Omega}$ is at most $2\delta$.

Having chosen $\psi$, we now choose $\phi$ with support so small that
$$
\lim_{\lambda \to \infty} \frac1{N(\lambda)}  \sum_{\lambda_j < \lambda} \langle \phi u_j^b, u_j^b \rangle_Y \leq (m^2 M)^{-1} \delta^2,
$$
where $m^2M$ is as in the proof of (ii) of Lemma~\ref{G}.
This is possible by \eqref{gaussian} and the Karamata Tauberian theorem.
It follows that there is a set of density $S''_\delta$ at least $1 - \delta$ where $\langle \phi u_j^b, u_j^b \rangle_Y \leq \delta$.

Now we write, from \eqref{uj-norm},
\begin{equation}
\| u_j^b \|_{L^2(Y)}^2 = \frac{\langle \psi^2 u_j, u_j \rangle_{\Omega}
- \langle S \phi u^b_j, u^b_j  \rangle_Y }{\langle S(1 - \phi) \hat u^b_j, \hat u^b_j  \rangle_Y }
\label{norm-eqn}\end{equation}
The denominator converges to a positive limit along a density one subsequence $S$, which by \eqref{inc} increases as either $\psi$ or $1 - \phi$ increases. On the other hand, along $S'_\delta \cap S''_\delta$, the numerator has oscillation at most $3\delta$. Hence along the subsequence $S \cap S'_\delta \cap S''_\delta$, the sequence of norms $\| u_j^b \|_{L^2(Y)}$ has oscillation at most a constant, depending only on $\Omega$, times $\delta$.

By Lemma~\ref{lem-seq}, this proves that the sequence of norms has a limit along a density one subsequence. To compute the limit, we let the support of $\psi$ increase to $\Omega$, and the support of $\phi$ shrink to $\Sigma$. Then $\langle \psi^2 u_j, u_j \rangle_{\Omega} \to 1$, $\langle S \phi u_j^b, u_j^b \rangle \to 0$, and by \eqref{step1}, \eqref{S-neu} and \eqref{inc}, the denominator of \eqref{norm-eqn} converges to
$$
\frac{\omega_{\Neu}(S(1 - \phi))}{\omega_{\Neu}(\Id)} \to \frac1{\omega_{\Neu}(\Id)}
$$
along a density one subsequence. Hence, $\langle u_j^b, u_j^b \rangle$ converges to $\omega_{\Neu}(\Id)$ along this subsequence.
This proves \eqref{step2}, and completes the proof of Theorem~\ref{main} for Neumann boundary conditions.


\section{The Robin boundary condition}\label{Neumann}

The proof for Robin boundary conditions, $\pa_\nu u = \kappa u$, for a real function $\kappa$ on the boundary, is almost identical. In this case, Green's formula gives us
\begin{equation*}
u(x) = \int_Y \frac{\pa G_0(x,y;\lambda)}{\pa \nu} u(y) - G_0(x,y;\lambda) \frac{\pa u}{\pa \nu} \, d\sigma(y), \quad x \in \text{ interior } \Omega.
\end{equation*}
Using the boundary condition we find
\begin{equation*}
u(x) = \int_Y \Big( \frac{\pa G_0(x,y;\lambda)}{\pa \nu} - \kappa(y)G_0(x,y;\lambda) \Big) u(y) \, d\sigma(y), \quad x \in \text{ interior } \Omega.
\end{equation*}
Let $E_h$ be the operator of Proposition~\ref{EandF}. Letting $x$ approach the boundary as before, we find that
$$
u_j^b = F^\kappa_{h_j} u_j^b, \quad F^\kappa_h = F_h + E _h \circ \kappa.
$$
By Proposition~\ref{EandF}, the operator $E \circ \kappa$ has
the same structure as $F$, but is one order lower; the FIO part of $E \circ \kappa$ is order $-1$, instead of zero. Hence, the principal symbol of the FIO part of $F^\kappa$ is {\it identical} to the case of Neumann boundary conditions. Consequently, the result, and the proof, are the same as for Neumann boundary conditions.


\section{Dirichlet boundary condition}

In this section, we modify the Neumann proof so that it works for
Dirichlet boundary conditions.

\subsection{Local Weyl law} First, we prove the local Weyl law
\begin{equation}
\lim_{\lambda \to \infty} \frac1{N(\lambda)} \sum_{\lambda_j \leq \lambda} \lambda_j^{-2} \langle A_{h_j} u_j^b, u_j^b \rangle \to \omega_{\Dir}(A)\label{Weyl-D}\end{equation}
for boundary values of
Dirichlet eigenfunctions. Recall that in this case the boundary
values are given by $ u_j^b =
\partial_{\nu} u_j |_{Y}$.

We prove \eqref{Weyl-D} in the same way as in
the Neumann case.
Reductions (i) - (v) work without any
modifications. Only the calculation in (vi) in which we consider
even order differential operators is somewhat different. With $A_h
= h^{2k}P $ where $P$ is a positive differential operator of
order $2k$ as in the Neumann case, we consider the slightly
different
$$
\tilde e_A(t) = \sum_{j=1}^\infty e^{-t\lambda_j^2} \lambda_j^{-2}
\langle A_h u_j^b, u_j^b \rangle .
$$
We have
$$
\big( \frac{d}{dt} \big)^{k+1} \tilde e_A (t) = \sum_{j=1}^\infty
e^{-t\lambda_j^2} \langle P^* P u_j^b, u_j^b \rangle \; = \; tr
P \tilde E(t).
$$
Hence if we apply $P D_x D_{x'}$ to $e^{-t \Delta_{\Dir}}$, and restrict to $x = x' =0$, we obtain a kernel whose trace is $(d_t)^{k+1}\tilde e_A(t)$.

Using the expression \eqref{Seel} for the heat kernel, with terms $\gamma_0$, $\delta_0$ now given by \eqref{top-terms-D}, we obtain
$$
\Big( \frac{d}{dt} \Big)^{k+1} \tilde e_A(t) = 2 \int e^{-t ( \eta^2 + |
\eta |^2_h)} |p(y,\eta)|^2  \xi^2 \, dy \, d\xi \, d\eta \ +
O(t^{-n/2 - k - 1/2}).
$$
This is the same expression as \eqref{app-1} (with $k=0$), except that we have an extra factor of $\xi^2$ from the two normal derivatives in $t$. Following the calculation in the appendix,
this extra factor has the effect of changing $m$ to $m+1$, and introducing the factor $1 - |\zeta|^2$ in the integral over $|\zeta| \leq 1$, changing the Neumann measure into the Dirichlet measure. Hence we end up with
$$
 \Big( \frac{d}{dt} \Big)^{k+1} \tilde e_A(t) \sim t^{-(\frac{n}{2} + k + 1)} (\frac{n}{2} + k)(\frac{n}{2} + k-1)\dots
(\frac{n}{2}) \frac{4}{(4\pi)^{\frac{n}{2}} \vol(S^{n-1})} \int\limits_{B^* Y} \sigma(A) \, d\mu_{\Dir}
$$
where $\sim$ means up to an error which is $O(t^{-(n/2 + m +
1/2)})$ as $t \to 0$. We conclude that
$$
\tilde e_A(t) = (4\pi t)^{-n/2}  \frac{4}{\vol(S^{n-1})} \int_{B^* Y} \sigma(A) \, d\mu_{\Dir}
+ O(t^{-(n-1)/2}), \quad t \to 0.
$$
Then, the Karamata Tauberian theorem gives \eqref{Weyl-D}.

\subsection{Quantum ergodicity}
We now sketch the proof of quantum ergodicity in the case of Dirichlet boundary conditions, which is the statement that
$$
\lambda_j^{-2} \langle A_{h_j} u_j^b, u_j^b \rangle \to \omega_{\Dir}(A),
$$
along a density one subsequence of integers $j$.

A crucial change from the Neumann case is that the formula
$F_{h_j} u_j^b = u_j^b$ in the Neumann case gets replaced by
\begin{equation} \label{DBE} F^*_{h_j} u_j^b = - u_j^b \end{equation}
in the Dirichlet case, due to the fact that the normal derivative
falls on the opposite coordinate. Hence we now need to use the adjoint of Lemma~\ref{Egorov-N}, which for convenience we state here.
\begin{lem}\label{Egorov-D} Suppose $A_h$ is a pseudodifferential operator of order zero, with $ \WF(A) \subset \mathcal{R}^N$.
Then, one can express
$$
(F_h)^N A_h (F_h^*)^N = B_h + S_h,
$$
where $\| S_h \|_{L^2 \to L^2} \leq C h$, and $B_h$ is a pseudodifferential operator of order zero with symbol
\begin{equation}
 \sigma(B) = \begin{cases}
\gamma^{-1}(q) \gamma(\beta^{-N}(q)) a(\beta^{-N}(q)) , \quad q \in \mathcal{R}^{-N} \\
0 , \phantom{ \gamma(q)^{-1}\gamma(\beta(q)) b(\beta^{(N)}(q)) } \quad q \notin \mathcal{R}^{-N}. \end{cases}
\label{Egorov-D-formula}\end{equation}
\end{lem}

Hence the Egorov operator $T$ (see section~\ref{ENDO})
gets replaced by
$$
T^* f(q) = \frac{\gamma(\beta^{-1}(q))} {\gamma(q)}
f(\beta^{-1}(q).
$$

{\it Step 1.} The first task is to choose a pseudodifferential
operator $\tilde R$ whose symbol is approximately equal to the invariant
function  $c \gamma (q)^{-1} = c (1 - |\eta|^2)^{-1/2}$, where $c$
is as in \eqref{c}. Given  $\epsilon_1 > 0$, we may choose $\tilde R$ so that its symbol is supported in $B^* Y$, and such that
\begin{equation}
\int \big| \sigma(R) - c \gamma (q)^{-1} \big|^2 \gamma (q)
d\sigma(q) < \epsilon_1. \label{R-symbol}\end{equation}
We then prove the following analogue of equation \eqref{step1}:
\begin{equation}
\limsup_{\lambda \to \infty} \frac1{N(\lambda)} \sum_{\lambda_j < \lambda}  \lambda_j^{-2} \Big( \langle Au_j^b, \hat u_j^b \rangle - \frac{\omega_{\Dir}(A)}{\omega_{\Dir}(\Id)} \| u_j^b \|_{L^2(Y)} \Big)^2 < C\epsilon,
\label{step1-D}\end{equation}
where $C$ depends only on $\Omega$ and $A$.
To do this, we define
\begin{equation}
{a}^{(N)}(q) = \frac1{N} \sum_{k=1}^N \frac1{\gamma(q)} { a(\beta^{-k}(q)) }{\gamma(\beta^{-k}(q))},
\label{aN-D}\end{equation}
an averaged version of the symbol in \eqref{Egorov-D-formula}.
(Here it is not necessary to introduce the function $a''$ as in the Neumann proof.) By the mean ergodic theorem,
$$
\gamma a^{(N)} \to c \omega_{\Dir}(A) \text{ in } L^2.
$$
In fact, since $\gamma a^{(N)}$ is uniformly bounded in $L^\infty$, we also get convergence in $L^p$ for $2 \leq p < \infty$, in particular for $p=6$. Since $\gamma^{-1/2} \in L^{3}(B^* Y)$, this implies that
$$
\gamma^{1/2} a^{(N)} \to c \omega_{\Dir}(A) \gamma^{-1/2} \text{ in } L^2.
$$
Thus, we can choose an integer $N$ such that
$$
\big\| \big( a^{(N)} \gamma - c \omega_{\Dir}(A)  \big) \gamma^{-1/2} \big\| < \epsilon.
$$
Choosing $a'$ as in section~\ref{Proof}, we obtain the analogue of \eqref{a'N-bound}:
\begin{equation}
\big\| \gamma^{1/2} \big( {a'}^{(N)} - c \gamma^{-1} \omega_{\Dir}(A) \big) \|_2^2 < 4\epsilon.
\label{a'N-bound-D}\end{equation}
The rest of Step 1 then follows the argument in section~\ref{Proof}.
\

{\it Step 2.}  Again the argument in section~\ref{Proof} adapts with minor changes. Since the boundary value $u_j^b$ here is the normal derivative of the eigenfunction, the value at an interior point is given by
$$
u_j(z) = \int_Y E_{h_j}(z,y) u_j^b(y) d\sigma(y),
$$
where $E_h$ is the free Green function at energy $h^{-2}$, instead of \eqref{Green}. Thus,
$$
\langle \psi^2 u_j, u_j \rangle = \langle E^*_{h_j} \psi^2 E_{h_j} u_j^b, u_j^b \rangle,
$$
so we need to analyze the operator $E^*_{h_j} \psi^2 E_{h_j} $. Lemma~\ref{G} gets replaced by

\begin{lem}\label{EE} Let $\phi$ be a smooth function on $Y$ which is identically one near $\Sigma$. Then

(i) $h^{-2} E_{h}^* \psi^2 E_{h} (1 - \phi)$ is a semiclassical pseudodifferential operator of order zero on $Y$;  and

(ii) For each $1 \leq i \leq m$, $h^{-2} E_{h}^* \psi^2 E_{h} $ is bounded on $L^2(Y)$, uniformly in $h$.

(iii) The limit of $\omega_{\Dir}(h^{-2} E^* \psi^2 E (1 - \phi))$, as the support of $\phi$ shrinks to $\Sigma$, is
\begin{equation}
\frac{1}{\vol(\Omega)} \int_\Omega \psi^2(z) \, dz.
\label{S-dir}\end{equation}
\end{lem}

The proof follows the same lines as the proof of Lemma~\ref{G}. In the proof of (iii), the lack of two normal derivatives in the operator is exactly compensated by the change in measure from $\mu_{\Neu} = \gamma^{-1} d\sigma$ to $\mu_{\Dir} = \gamma d\sigma$.

The proof is completed by writing
\begin{equation}
\| u_j^b \|_{L^2(Y)}^2 = \frac{\langle \psi^2 u_j, u_j \rangle_{\Omega}
- \lambda_j^{-2} \langle S \phi u^b_j, u^b_j  \rangle_Y }{ \lambda_j^{-2} \langle S(1 - \phi) \hat u^b_j, \hat u^b_j  \rangle_Y }
\label{norm-eqn-D}\end{equation}
and using Lemma~\ref{lem-seq} and the result of Step 1, as in section~\ref{Proof}.


\section{$\Psi^1$-Robin boundary condition}\label{Psi}

In this section we consider boundary conditions of the form
$$
\pa_{\nu} u = K u,
$$
where $K$ is a self-adjoint pseudodifferential operator of order $1$ with nonnegative principal symbol on the boundary. We shall denote the operator $H_K$ in this case. Recall that it is defined by a quadratic form \eqref{quadform}.

The local Weyl law for $H_K$ is as follows.

\begin{lem} \label{LWLN} Suppose we impose the boundary condition
above. Let $k$ denote the principal symbol of $K$.
 Let $A_h$ be a semiclassical operator of order
 zero on $Y$. Then
\begin{equation}
\lim_{\lambda \to \infty} \frac{1}{N(\lambda)} \sum_{\lambda_j \leq \lambda}
\langle A_{h_j} u^b_j, u^b_j \rangle =  \omega_K(A) \equiv \frac{2}{\vol(S^{n-1})\vol(\Omega)} \int\limits_{B^* Y} \frac{ a(q) \gamma(q) }{ \gamma(q)^2 + k(q)^2 } \, d\sigma(q).
\label{Weyl-limit-K}\end{equation}
\end{lem}

\begin{proof} Identical in structure to the proof of Lemma~\ref{LWL}.
\end{proof}

Notice that if $K = 0$ then this reproduces the
result of the previous section. Consider now (for a manifold with smooth boundary) $K = \alpha \Delta_Y^{1/2}$.
As $\alpha \to \infty$, this boundary condition approaches Dirichlet, in some sense. To compare to our Dirichlet result, we need to account for the fact that
here we are considering the $L^2$ norm of the function at the boundary,
instead of the normal derivative. To change to the normal derivative is
equivalent, by our boundary condition, to changing to $\alpha \Delta_Y^{1/2} u$,
so to compare, we need to multiply the weight factor in \eqref{Weyl-limit-K} by $\alpha^2 |\eta|^2$. If we do this, then as $\alpha \to \infty$ we get
$$ \lim_{\alpha \to \infty}
\int_{B^* Y} \frac{ a(q) \gamma(q) \alpha^2 |\eta|^2 }{ \gamma(q)^2 + \alpha^2 |\eta|^2 } \, d\sigma(q)
= \int_{B^* Y} a(q) \gamma(q)  \, d\sigma(q),
$$
which is the Dirichlet measure. Thus, for this family of boundary conditions
the corresponding measures $\omega_{\alpha K}$ interpolate between the Neumann and Dirichlet measures.

In order to adapt the arguments above to the $\Psi^1$-Robin boundary conditions, we need an operator, which is essentially a FIO, which leaves the $u_j^b$ invariant.

\begin{lem}
Let $\epsilon > 0$ and $N$ be given. Then there is a subset $V$ of $B^*Y \times B^*Y$ such that both projections $V_1$ and $V_2$ onto the left and right factors of $B^*Y$ have measure less than $\epsilon$, and a operator $F^K_h$, acting on $L^2(Y)$, depending on the parameter $h > 0$, such
 that

(i) $F^K$ leaves the $u_j^b$ invariant:
$$
F^K_{h_j} u_j^b = u_j^b,
$$

(ii) there is a decomposition
$$
F^K = F^K_1 + F^K_2 + F^K_3
$$
such that $(F^K_1)_h$ is a semiclassical FIO of order $0$, associated to the billiard map, $(F^K_2)_h$ is a semiclassical pseudodifferential operator of order $-1$, and $(F^K_3)_h$ has wavefront set contained in $V$;

(iii) the symbol of $F^K_1$ in $\Cbill \setminus V$ is given by
$$
 \sqrt{\frac{\gamma(\beta(q))}{\gamma(q)}}  \frac{ \gamma(q) - i k(q) }{ \gamma(\beta(q)) + i k(\beta(q)) };
$$

(iv) If $A_h$ is a semiclassical pseudodifferential operator with microsupport disjoint from $V_2$, then for any $-N \leq m \leq N$, the operator $((F^K)^*)^m A (F^K)^m = A^{(m)} + S$, where $A^{(m)}$ is a semiclassical pseudodifferential operator with symbol
$$
\frac{ \gamma(\beta^m(q)) } { \gamma(\beta^m(q))^2 + k(\beta^m(q))^2 } \frac{\gamma(q)^2 + k(q)^2 }{ \gamma(q)  } a(\beta^m(q))
$$
and $\| S_h \|_{L^2(Y) \to L^2(Y)} = O(h)$.

\end{lem}

\begin{proof} We first prove these assertions in the case $N=1$. We start from Green's formula:
$$
u(x) = \int_Y \frac{\pa G_0(x,y;\lambda)}{\pa \nu} u(y) - G_0(x,y;\lambda) \frac{\pa u}{\pa \nu} \, d\sigma(y), \quad x \in \text{ interior } \Omega.
$$
Using our boundary condition, we find
$$
u(x) = \int_Y \frac{\pa G_0(x,y;\lambda)}{\pa \nu} u(y) - G_0(x,y;\lambda) (Ku) \, d\sigma(y).
$$
Now we let $x$ tend to the boundary. From the normal derivative of $G_0$ we get $\half \delta_x(y)$ plus half the kernel $F_h$, $h = \lambda^{-1}$.
Using Proposition~\ref{F-corners}, given any neighbourhood $U$ of the set
$\Xi \times (\mathcal{R}^1)' \cup (\mathcal{R}^{-1})' \times \Xi$, we can decompose $F = F_1 + F_2 + F_3$, where $F_1$ is an FIO, $F_2$ is a pseudodifferential operator of order $-1$ and $F_3$ has operator wavefront set contained in $U$.

The other term is more complicated, since we are composing the `single layer potential' $E$, given by the kernel of $2G_0$ restricted to the boundary in both variables, with $K$. It is necessary to introduce a more elaborate decomposition of $E$ in this case, since $K$ is not a semiclassical operator, and there is no nice calculus containing both homogeneous and semiclassical operators. However, if $T$ is a semiclassical FIO whose operator wavefront set is disjoint  from the zero section in the right variable, that is, if
$$
\WF(T) \subset \{ (q, q') \mid q' \neq 0 \},
$$
then $TK$ is a semiclassical FIO with symbol $t(q, q') k(q')$. Here, $k$ is the {\it principal} symbol of $K$. This follows from Taylor \cite{T}, chapter 8, section 7, for example.

Thus, we let $V$ be any neighbourhood of $U \cup (Z  \cup \beta(Z)) \times Z$, where $Z$ is the zero section of $B^*Y$. (Notice that both $U$ and $V$, as well as both their projections, may have arbitrarily small measure.) We decompose $E = E_1 + E_2 + E_3$, where $E_1$ is an FIO of order $-1$ microsupported away from the zero section in the right variable, $E_2$ is a pseudodifferential operator of order $-1$ microsupported away from the zero section in the right variable, and $\WF(E_3) \subset V$.

Thus, we have
\begin{equation}
u_j^b = \int_Y \big( F_1 + F_2 + F_3 \big) u_j^b + \big( E_1 + E_2 + E_3 \big) K u_j^b .
\end{equation}
By Lemma~\ref{EandF}, the symbol of $E_2$ inside $B^*Y$, and outside $V_2 = \operatorname{proj}_2(V)$, is $-i \gamma^{-1}$. Thus, $E_2 K$ is a pseudodifferential operator of order $0$ with principal symbol $-i \gamma^{-1} k$ outside $V_2$. We rewrite this equation as
\begin{equation}
\Big( \Id - E_2 K \Big) u_j^b = \int_Y \big( F_1 + F_2 + F_3 \big) u_j^b + \big( E_1 + E_3 \big) K u_j^b .
\end{equation}
Since the symbol of $E_2 K$ is imaginary, there is a pseudodifferential operator $T$ such that $T (\Id - E_2 K) = \Id - E_4$, where $E_4$ is order $-1$. Moreover, the symbol of $T$ on
$B^*Y \setminus V_2$ is $(1 + i\gamma^{-1} k)^{-1}$. Thus, we get
\begin{equation}
 u_j^b = T \big( F_1 + F_2 + F_3 \big) u_j^b + T \big( E_1 + E_3 \big) K u_j^b  + E_4 u_j^b.
\end{equation}
Now we let
\begin{equation}\begin{aligned}
F^K_{1}(x,y) &= T F_1 + T E_1 K \\
F^K_{2}(x,y) &= T F_2 + E_4 \qquad\text{ and } \\
F^K_{3}(x,y) &= T F_3 + T E_3.\end{aligned}\end{equation}
Then we have
\begin{equation}
u^b_j = (F^K_{1} + F^K_{2} + F^K_{3})_{h_j} u_j^b, \quad h_j = \lambda_j^{-1},
\end{equation}
where $F^K_1$ is a semiclassical FIO with principal symbol
$$
\frac1{1 + ik(\beta(q))/\gamma(\beta(q))} \bigg( \! \sqrt{\frac{\gamma(q)}{\gamma(\beta(q))}} - \frac{  i k(q)}{\sqrt{\gamma(q) \gamma(\beta(q))}} \! \bigg)
= \sqrt{\frac{\gamma(\beta(q))}{\gamma(q)}}  \frac{ \gamma(q) - i k(q) }{ \gamma(\beta(q)) + i k(\beta(q)) }
$$
on $\Cbill \setminus V$;
$F^K_2$ is a semiclassical pseudodifferential operator of order $-1$, supported outside $V_2$; and $F^K_3$ is an error term with microsupport contained in $V$. This proves (i), (ii) and (iii).
Consequently, the symbol of $(F^K_1)^* A F^K_1$, for $A$ microsupported within $V_2$ is
$$
\frac{ \gamma(\beta(q)) } { \gamma(\beta(q))^2 + k(\beta(q))^2 } \frac{\gamma(q)^2 + k(q)^2 }{ \gamma(q)  } a(\beta(q)).
$$
Statement (iv), for $N=1$, then follows from the calculus of wavefront sets and the symbol calculus for FIOs, as in Proposition~\ref{Eg}.

To prove statement (iv) for arbitrary $N$ we need to enlarge the set $V$ to a set
\begin{equation*}
V^{(N)} = \{ (q, q') \mid \ \exists \ q'' \text{ such that } (q'', q') \in V, q = \beta^m(q'') \text{ for some } 0 \leq m \leq N-1 \}.
\end{equation*}
The projections onto the first and second factors of this set have small measure if $V$ has projections of small measure, since $\beta$ is measure-preserving. Then (iv) follows as in the proof of Proposition~\ref{Eg}.
\end{proof}

This Lemma allows us to run Step 1 of the argument in section~\ref{Proof} for the $\Psi^1$-Robin boundary condition.

To run Step 2, we again run into complications caused by composition of homogeneous and semiclassical pseudodifferential operators.  In the present case, we get instead of \eqref{interior}
\begin{equation}
\langle \psi^2 u_j,  u_j \rangle_{\Omega} =
 \langle (G_{h_j} + E_{h_j} K)^* \psi^2 (G_{h_j} + E_{h_j} K) u_j^b, u_j^b \rangle_Y .
\label{interior-K}\end{equation}

Here we need to analyze the operator $(G_{h} + E_{h} K)^* \psi^2 (G_{h} + E_{h} K) $.

\begin{lem}\label{GEK} Let $\phi$ be a smooth function on $Y$ which is identically one near $\Sigma$, and let $Q$ be a semiclassical pseudodifferential operator with microsupport disjoint from the zero section. Then

(i) $(G + E Q K)^* \psi^2 (G + E Q K) (1 - \phi)$ is a semiclassical pseudodifferential operator of order zero on $Y$;  and

(ii) For each $1 \leq i \leq m$, $(G + E Q K)^* \psi^2 (G + E Q K)$ is bounded on $L^2(Y)$, uniformly in $h$.

(iii) The limit of $\omega_{\Neu}((G + E Q K)^* \psi^2 (G + E Q K)(1 - \phi))$, as the support of $\phi$ shrinks to $\Sigma$, and the symbol of $Q$ increases to $1$, is
\begin{equation}
\frac{1}{\vol(\Omega)} \int_\Omega \psi^2(z) \, dz.
\label{S-neu-K}\end{equation}
\end{lem}

The proof of this Lemma follows the same lines as that of Lemma~\ref{G} and Lemma~\ref{EE}. The novelty here is in dealing with the cutoff operator $Q$, which gives us additional error terms
$$\begin{gathered}
\langle u_j^b, G^* \psi^2 E(\Id - Q)K u_j^b \rangle,
\langle u_j^b, ( E (\Id - Q)K )^* \psi^2 (E Q K) u_j^b \rangle,\\
\text{ and }
\langle u_j^b, (E(\Id - Q)K )^* \psi^2 (E(\Id - Q)K ) u_j^b \rangle
\end{gathered}$$
and their adjoints. We need to show, with the first term for example, that
$$
\lim_{\lambda \to \infty} \frac1{N(\lambda)} \sum_{\lambda_j < \lambda} \langle u_j^b, G^* \psi^2 E (\Id-Q) K u_j^b \rangle = 0.
$$
Writing this as
$$\begin{gathered}
\limsup_{\lambda \to \infty} \frac1{N(\lambda)} \sum_{\lambda_j < \lambda} \langle (\Id-Q) E^* \psi^2 G u_j^b,  K u_j^b \rangle  \\
\leq \limsup_{\lambda \to \infty} \frac1{N(\lambda)} \sum_{\lambda_j < \lambda} \epsilon^{-1} \| h^{-1} (\Id-Q) E^* \psi^2 G u_j^b \|^2 + \epsilon \|  h K u_j^b \|^2
\end{gathered}$$
we see that this can be made arbitrarily small by choosing $Q$ suitably, since the first term is controlled by $\omega_K((\Id-Q)(h^{-1} E^*) G)$ which can be made arbitrarily small, and the second can be made small by choosing $\epsilon$ small. To deal with the third term, we need to commute $K$ and $(E(\Id - Q))^* \psi^2 E (\Id - Q)$. This again causes difficulties since one is a homogeneous, and one a semiclassical, operator. However,  Lemma~\ref{commuting} allows us to commute the operators up to an error with an $o(1)$ operator norm as $h \to 0$. This allows us to complete Step 2 for the $\Psi^1$-Robin boundary condition, which completes the proof of quantum ergodicity in this case.


\section{Nonconvex domains}\label{nonconvex}
Here we adapt the argument to nonconvex domains. The problem with the argument above for nonconvex domains is that  the FIOs $E_{1,h}$ and $F_{1,h}$ from Section~\ref{Structure} have a canonical relation which is larger than the billiard relation, since it relates points on the boundary which are connected by a straight line even if the line passes outside the domain. In this section, we show that one can manufacture an invariant operator $F_h$ which satisfies the conditions of Proposition~\ref{F-corners}. The argument then proceeds as before.

We do this by modifying the metric outside the domain, while keeping it Euclidean inside. Let $b$ be a smooth, compactly supported nonnegative function on $\RR^n$ which vanishes on $\Omegabar$. Consider the metric
$$
g_s = (1 + sb) g_{\operatorname{Euclidean}} \ \text{ on } \RR^n.
$$
For sufficiently small $s$, no geodesics of $g_s$ starting at a point in $\Omegabar$ have conjugate points in $\Omegabar$. Let $G_s(z, z, h) = (\Delta_s - (h^{-1} + i0)^2)^{-1}(z, z')$ denote the kernel of the outgoing resolvent of the Laplacian on $\RR^n$ with respect to the metric $g_s$.

\begin{prop}\label{G-str} The kernel $G_t(z, z', h)$ may be represented
$$
G_t(z, z', h) = h^{2-n} e^{i\operatorname{dist}_s(z, z')/h} a(z, z', h), \quad z \neq z'.
$$
\end{prop}

\begin{proof}
We obtain the resolvent from the heat kernel by integrating against the wave kernel
$$
(\Delta_s - (\lambda + i0)^2)^{-1} = \int_{_\infty}^\infty e^{it \sqrt{\Delta_s}} e^{i |t| \lambda} \, dt.
$$
Recall that the wave kernel $e^{it \sqrt{\Delta_s}}$ is a family of FIOs, with
canonical relation given by
\begin{multline*}
\{ (z, \zeta, z', \zeta') \mid \ \exists \ \text{ a $g_s$-geodesic $\gamma$ of length $t$  with }  \\ z = \gamma(0) , \zeta = c \dot \gamma(0), z' = \gamma(s), \zeta' = -c\dot \gamma(s) \}.
\end{multline*}
It can therefore be written in the form
$$
e^{it \sqrt{\Delta_s}}(z, z')  = \int e^{i((z - Z'(z, z',s)) \cdot \zeta - t |\zeta|)} a(t, z, z', \zeta) \, d\zeta,
$$
where $a$ is smooth in all variables and symbolic of order $-1$ in $\zeta$,
and, $Z'(z, z',s)$ is the coordinate of $z'$ in $g_s$-normal coordinates centred around $z$.  Putting $\zeta = \lambda \xi$, we find that
$$
(\Delta_s - (\lambda + i0)^2)^{-1} =\lambda^{n-2}  \int_{-\infty}^\infty dt \int   d\xi \,  e^{i \lambda ((z - Z'(z, z',s)) \cdot \xi - t |\xi|)} e^{i |t| \lambda} a(t, z, z', \xi, \lambda^{-1})  .
$$
For $|t| \geq T_0$, smoothness of the wave kernel on $\Omegabar \times \Omegabar$ implies, by repeated integrations by parts in $t$, that the corresponding kernel is residual. For $0 < t_0 \leq |t| \leq T_0$, where $T_0$ is larger than the $g_s$-diameter of $\Omegabar$ for each $s$, the integral is a semiclassical FIO with canonical relation \begin{equation*}
\{ (z, \zeta, z', \zeta') \mid \ \exists \ \text{  $g_s$-geodesic $\gamma$  with }   z = \gamma(0) , \zeta =  \dot \gamma(0), z' = \gamma(s), \zeta' = -\dot \gamma(s) \}.
\end{equation*}
This is parametrized by the global phase function $\operatorname{dist}_s(z, z')$.
Finally, for $|t| < t_0$, where $t_0$ is sufficiently small, the wave kernel $e^{it \sqrt{\Delta_s}}$ is identical with the free wave kernel $e^{it \sqrt{\Delta_0}}$ due to finite propagation speed, so $G_s(z, z', h)$ is microlocally identical with the free resolvent kernel near the diagonal.   This proves the lemma.
\end{proof}

To deal with nonconvex domains, we define
$$
F^s_h(y, y') = 2 \dbyd{}{\nu_y} G_s(y, y', h), \quad y, y' \in Y,
$$
and average over $s$ to obtain the operator
\begin{equation}
\tilde F_h = \int_0^1 \chi(s) F^s_h \, ds.
\label{average-op}\end{equation}
Here $\chi$ is smooth, supported in $
(0, \delta)$, and nonnegative with integral $1$. The averaged operator $\tilde E_h$ is defined similarly.

Let us say that a point $(y, \eta, y', \eta')$, with $\eta = d_y |y - y'|$, $\eta' = -d_{y'} |y - y'|$ is a \emph{spurious} point of $\WF(F_h)$ if the line $\overline{yy'}$ leaves $\Omegabar$. We also recall the definition of $\Xi$ in \eqref{Xi}. Then the following analogue of Proposition~\ref{F-corners} holds.

\begin{prop}\label{averaging}
(i) The operator $\tilde F_h$ is an invariant operator for the Neumann eigenfunctions:
$$
\tilde F_h u_j^b = u_j^b.
$$
Moreover, the analogous statements for other boundary conditions are valid  for the corresponding averaged invariant operators.

(ii) Let $U$ be any neighbourhood of $\Xi \times (\mathcal{R}^{1})' \, \cup \, (\mathcal{R}^{-1})' \times \Xi \, \cup \, (\mathcal{R}^{-1})' \times (\mathcal{R}^{1})' $.  Then for a suitable averaged invariant operator $\tilde F_h$, there is a decomposition
$$
\tilde F_h = \tilde F_{1,h} + \tilde F_{2,h} + \tilde  F_{3,h},
$$
where $\tilde F_1$ is a Fourier Integral operator of order zero associated with the canonical relation $\Cbill$, $\tilde F_2$ is a pseudodifferential operator of order $-1$ and $\tilde F_3$ has operator wavefront set contained in $U$.

(iii)  The symbol of $\tilde F_{1,h}$ is the same as for $F_h$.
\end{prop}

\begin{proof}

(i) Each $F^t_h$ is an invariant operator since the metrics $g_t$ all coincide with the Euclidean metric inside $\Omega$. Hence any average of the $F^t_h$ is also an invariant operator.

(ii) If $b > 0$ somewhere on $\overline{y y'}$, then
$d/ds  (\operatorname{dist}_s(y, y')) > 0$ for small $s$
by the positivity of $b$. Hence, in the integral \eqref{average-op}, the phase is locally non-stationary, implying that the kernel is locally residual.

A compactness argument shows that if the support of $b$ is sufficiently close to $\Omega$, then any spurious vector $q \in \WF(F_h)$ is either in $U$ or the line $\overline{yy'}$ joining the endpoints meets the support of $b$, in which case $q \notin \WF(\tilde F_h)$ by the remarks of the previous paragraph. 

As for the non-spurious wavefront set, the decomposition follows from the structure of $G_s(z, z', h)$ given in Proposition~\ref{G-str} and the proof of Proposition~\ref{F-corners}.

(iii) Let $(y, \eta, y', \eta') \in \WF(F_{h})$, such that the line joining $y, y'$ remains inside $\Omega$. Then the metrics $g_s$ are identical in a neighbourhood of this line, hence
the wave kernels $e^{it \sqrt{\Delta_s}}$, and therefore,
the resolvent kernels $G_s$ are all identical modulo residual terms near $(y, y')$. Hence all the $F^s_h$ are microlocally identical near $(y, \eta, y', \eta')$, which implies (iii).

\end{proof}

Lemma~\ref{Egorov-N} now holds for a suitable averaged invariant operator $\tilde F_h$. Hence, the arguments of Section~\ref{Proof} go through, proving the result in the nonconvex situation.


\section{Appendix. The heat kernel}

We begin by exhibiting the asymptotic expansion of the heat kernel for a domain with smooth boundary, first under $\Psi^1$-Robin boundary conditions and then under Dirichlet boundary conditions.

\subsection{Neumann or $\Psi^1$-Robin boundary condition}

The asymptotic expansion of the heat kernel as $t \to 0$ may be obtained from the asymptotic expansion of the resolvent $(\Delta_B - \lambda)^{-1}$ as $\lambda \to \infty$ within a sector disjoint from the positive real axis. The construction of Seeley in \cite{Seeley}, \cite{Seeley1} goes through for a boundary condition of the form
\begin{equation}
\pa_{\nu} u = Ku, \quad K \in \Psi^1(Y), \quad k = \sigma_1(K) \geq 0.
\label{K-bc2}\end{equation}
In local coordinates, the Seeley parametrix takes the form
\begin{equation}\begin{gathered}
R_N(\lambda) = \sum_{j=0}^{N} \Bigg\{
\int e^{i(x-x')\cdot \xi} e^{i(y - y') \cdot \eta} c_{-2-j}(x, y, \xi, \eta, \lambda) \, d\xi   \, d\eta +  \\
 \int e^{-i x' \xi} e^{i(y-y')\cdot \eta}  d_{-2-j}(x, y, \xi, \eta, \lambda) \,  d\xi \, d\eta \Bigg\}.
\end{gathered}
\end{equation}
We obtain the heat kernel $e^{-t\Delta_B}$
by performing a contour integral:
$$
e^{-t \Delta_B} = \frac{i}{2\pi} \int_\Gamma e^{-t\lambda} (\Delta_B - \lambda)^{-1} d\lambda,
$$
where $\Gamma$ `encloses' the spectrum of $\Delta_B$, eg
$$
\Gamma = \{ -L + s e^{-i\theta}, s \in [0, \infty) \} \cup \{ -L + s e^{i\theta}, s \in [0, \infty) \}, \quad 0 < \theta < \pi.
$$
Substituting the series for $R_N(\lambda)$ in for the resolvent, we obtain a parametrix for the heat kernel of the form in a single coordinate patch
\begin{equation}\begin{gathered}
\sum_{j=0}^{N} \Bigg\{
(2\pi)^{-n} \int e^{i(x-x')\cdot \xi} e^{i(y-y') \cdot \eta} \gamma_j(x, y, \xi, \eta, t) \, d\xi  \, d\eta \ + \\
(2\pi)^{-n} \int e^{i(y-y')\cdot \eta} e^{-i x' \xi} \delta_j(x, y, \xi, \eta, t) \, d\xi \, d\eta \Bigg\}.
\end{gathered}\label{Seel}
\end{equation}
After patching and summing over coordinate patches, we denote this operator $Z_N(t)$.
This is a good approximation to the heat kernel as $t \to 0$, in the sense that
the kernel of the difference,
\begin{equation}
q_N(t,x,y) = e^{-t \Delta_B}(x,y) -  Z_N(t)(x,y),
\label{heat-par}\end{equation}
satisfies estimates
\begin{equation}
| D_x^\alpha D_y^\beta q_N(t) | \leq C t^{-(n + \alpha + \beta) + N + 1)/2} e^{-\delta |x-y|^2/t}
\label{heat-par-est}\end{equation}
for some $\delta > 0$.
We will be interested mostly in the top terms in the expansion of the heat kernel, corresponding to $\gamma_0$ and $\delta_0$. To write them in the simplest possible way, we assume that we have chosen coordinates so that
$|\pa_{x}| = 1$ and $\pa_{x}$ and $\pa_{y_i}$ are orthogonal at the boundary, and so that the lines $\{ y = \text{constant} \}$ are geodesics close to the boundary. Then, the metric $g^{ij}$ takes the form
$$
g^{ij} = \begin{pmatrix}
\phantom{gg} & & \phantom{gg}& 0 \\
& h^{ij} & & 0 \\
& & & \cdot \\
\ \dots & 0 & \dots & 1
\end{pmatrix}
$$
where $h$ is the induced metric on the boundary.
Let us denote $\sum_{i,j \leq n-1} h^{ij} \eta_i \eta_j$ by $|\eta|^2_h$. Then
the top two coefficients $c_2$ and $d_2$ in the parametrix for the resolvent with boundary conditions \eqref{K-bc2} are
\begin{equation}\begin{aligned}
c_2(x,y,\xi, \eta, \lambda) &= (\xi^2 + |\eta|_h^2 - \lambda)^{-1} \\
d_2(x,y,\xi, \eta, \lambda) &= \frac{i\xi - k(y,\eta)}{\sqrt{|\eta|_{h_0}^2 - \lambda} + k} \ \frac{e^{-\sqrt{|\eta|_{h_0}^2 - \lambda}\ x}}{\xi^2 + |\eta|_{h_0}^2 - \lambda}, \quad \operatorname{Re} \sqrt{|\eta|_{h_0}^2 - \lambda} > 0.
\end{aligned}\end{equation}
Recall that $k(y, \eta)$ is the {\it principal} symbol of $K$, hence is nonnegative by assumption, and homogeneous of degree one. Hence
$d_2$ has no more singularities than it does for the Neumann boundary condition $k \equiv 0$. Essentially for this reason, the Seeley parametrix works just as well for $\Psi^1$-Robin boundary conditions as it does for  Neumann boundary condition.

The corresponding top terms $\gamma_0$ and $\delta_0$ for the heat parametrix at $x=0$ are
\begin{equation}\begin{aligned}
\gamma_0(0, y, \xi, \eta, t) &= e^{-t ( \xi^2 + | \eta |^2_h)} \\
\delta_0(0, y, \xi, \eta, t) &= e^{-t (\xi^2 + | \eta |^2_h)} \frac{i\xi + k(y,\eta)}{i\xi - k(y,\eta)}.
\end{aligned}\label{top-terms}\end{equation}

Next, for piecewise smooth domains, we show

\begin{lem}\label{HK} Suppose that $\Omega \subset \RR^n$ is a piecewise smooth domain.  Then for $\Psi^1$-Robin boundary conditions, the heat kernel $e^{-t\Delta_K}(z,z')$ admits Gaussian bounds for $t \leq 1$,
\begin{equation}
| e^{-t\Delta_K}(z,z') | \leq C t^{-n/2} e^{-\delta |z-z'|^2/t}
\label{gaussian}\end{equation}
for some constants $C, \delta$ depending on $\Omega$. Moreover, in the regular part $Y^o$ of the boundary, the local expansion is valid.
\end{lem}

\begin{proof}
The operator $\Delta_K$ is the self-adjoint operator associated with the quadratic form \eqref{quadform}. When $K$ is nonnegative, this quadratic form is larger than the corresponding form for the Neumann boundary condition ($K \equiv 0$), so  \eqref{gaussian} follows from \cite{D}, Theorems 2.4.4 and 3.2.9.

To prove the second part of the Lemma, we fix a point $y_0 \in Y^o$ and extend $\Omega$ to a manifold $\Omega' \supset \Omega$ with smooth boundary so that the boundaries of $\Omega$ and $\Omega'$ coincide on the support of the kernel of $K$ or $\kappa$, and on a neighbourhood of $y_0$.  We can then choose a
boundary condition $B$ for $\Omega'$ which agrees with the original boundary condition wherever the two boundaries coincide. Since $H_{\Omega'}$ has local asymptotics, it is enough to show that $H_{\Omega'} - H_\Omega$ has trivial asymptotics near $y_0$. 

The heat kernel $H_\Omega(t, z, z')$ for $\Omega$ for the boundary condition $B$ can be expressed in terms of that for $\Omega'$ by
$$
H_\Omega(t, z, z') = H_{\Omega'}(t, z, z') - \int_0^t \int_Y H_\Omega(t-s, z, y) \big( B H_{\Omega'}(x, y, z') \big) \, d\sigma(y) \, ds.
$$
Notice that $B H_{\Omega'}(x, y, z') = 0$ for $y$ near $y_0$. Hence the Gaussian bounds for $H_\Omega$ and $H_{\Omega'}$ (together with all derivatives of $H_{\Omega'}$) imply that the difference between the heat kernels has trivial asymptotics, as required. 
\end{proof}

\

Our next task is to show that for $\Psi^1$-Robin boundary conditions, the operator $E(t)$, obtained by restricting the heat kernel to the boundary in both variables, is trace class.

To do this, we prove the bound
\begin{equation}
\| u_j^b \|_{L^2(Y)} \leq C \lambda_j^2.
\label{psi-bound}\end{equation}
Note that $\Delta_K$ has no null space since $K$ is positive, and $\Delta_K$ has discrete spectrum since the quadratic form is larger than the Neumann quadratic form. Hence $\Delta_K^{-1}$ exists, and it is certainly a bounded map from $L^2(\Omega)$ into the form domain $H^1(\Omega)$. We conclude that the eigenfunction $u_j$, which satisfies
$$
u_j = \Delta_K^{-1} (\lambda_j^2 u_j),
$$
satisfies $\| u_j \|_{H^1(\Omega)} \leq C \lambda_j^2$. Hence by the Sobolev trace theorem, $\| u_j^b \|_{L^2(Y)} \leq C \lambda_j^2$.
Since the kernel of $E(t)$ is given by
\begin{equation}
E(t) = \sum_j e^{-t \lambda_j^2} u^b_j(y) u^b_j(y'),
\label{Ekernel}\end{equation}
and $\sum_j e^{-t \lambda_j^2} \lambda_j^2$ converges for each $t > 0$, it follows that $E(t)$ is trace class for all $t > 0$.

Moreover, if $A$ is a differential operator on $Y$ with coefficients supported inside $Y^o$, of order $k$, it is easy to show that
$$
\| Au_j^b \|_{L^2(Y)} \leq C \lambda_j^{k+2}.
$$
Since $\sum_j e^{-t \lambda_j^2} \lambda_j^{k+1}$ converges for any $k$, we have by Lidskii's theorem that $A E(t)$ is trace class, and the trace is given by the integral of the kernel of $A E(t)$ on the diagonal.

\vskip 5pt

Next, we calculate the trace of $P E(t)$, where $P$ is a positive, order $2k$ differential operator on $Y$ with symbol $p(y, \eta)$. Decomposing $e^{-t H_K} = H_0(t) + q_0(t)$ as in \eqref{heat-par}, we first substitute for $E(t)$ the operator obtained by restricting the kernel of $H_0$ to the boundary in both variables. This we can compute directly: we get
\begin{equation}
(2\pi)^{-n} \int e^{-t ( \xi^2 + | \eta |^2_h)} p(y, \eta) \frac{2i\xi}{i\xi - k(y,\eta)} \, dy \, d\xi \, d\eta \ .
\label{app-1}\end{equation}
To compute this we make a linear change of variable from $\eta$ to $\zeta = \zeta(y,\eta)$ which is an orthonormal basis for $h^{ij}(y)$ for every $y$. Write $\mu = \det (h^{ij})^{-1/2}$, so that the Riemannian measure on $Y$ is $\mu \, dy$. Under this transformation, $p(y,\eta)$ changes to $\tilde p_{2k}(y, \zeta) + \tilde p_{2k-1}(y, \zeta)$, say, where $\tilde p_{2k}$ is homogeneous of degree $2k$ in $\zeta$ and $\tilde p_{2k-1}$ is a polynomial of degree $2k-1$ in $\zeta$.  We obtain
\begin{equation*}\begin{gathered}
(2\pi)^{-n} \int e^{-t ( \xi^2 + | \zeta|^2)} \Big( \tilde p_{2k}(y, \zeta) + \tilde p_{2k-1}(y, \zeta) \Big) \frac{2i\xi}{i\xi - k}  \mu \, dy \, d\xi \, d\zeta \\
= (2\pi)^{-n} \int e^{-t ( \xi^2 + | \zeta|^2)} \Big( \tilde p_{2k}(y, \zeta) + \tilde p_{2k-1}(y, \zeta) \Big) \frac{\xi^2 - 2i\xi k}{\xi^2 + k^2}  \mu \, dy \, d\xi \, d\zeta .
\end{gathered}\end{equation*}
The $-2i\xi k$ term can be dropped because this gives an odd integral in $\xi$. If we then change to polar coordinates, $r^2 = \xi^2 + |\zeta|^2$, then we get
\begin{equation}
(2\pi)^{-n} \int_0^\infty e^{-t r^2} r^{n-1+2k} \, dr \int \mu \, dy \int\limits_{|(\xi, \zeta)| = 1} \tilde p_{2k}(y, \zeta) \frac{2\xi^2}{\xi^2 + k^2} + \dots
\end{equation}
Here the dots represent the contribution of $\tilde p_{2k-1}$. If we integrate by parts $k$ times, we find that
$$
 \int_0^\infty e^{-t r^2} r^{n-1+2k} \, dr  =  (n/2 + k-1) (n/2 + k-2) \dots (n/2)t^{-k} \int_0^\infty e^{-t r^2} r^{n-1} \, dr .
$$
Turning this back into an integral over $\RR^n$ gives
$$
\int_0^\infty e^{-t r^2} r^{n-1} \, dr = \frac1{\vol(S^{n-1})} \int_{\RR^n} e^{-t |x|^2} \, dx = \frac1{\vol(S^{n-1})} \pi^{n/2} t^{-n/2}.
$$
Writing $\xi^2 = 1 - |\zeta|^2$ on the unit sphere and representing the measure as $(1 - |\zeta|^2)^{-1/2} d\zeta$, we get
\begin{equation*}
\frac{(\frac{n}{2} + k-1)(\frac{n}{2} + k-2)\dots (\frac{n}{2}) t^{-n/2-k}}{(4\pi)^{n/2} \vol(S^{n-1})} \int \mu \, dy \int\limits_{|\zeta| \leq 1} \tilde p_{2k}(y, \zeta) \frac{4\sqrt{1 - |\zeta|^2}}{1 - |\zeta|^2 + k^2} \, d\zeta + \dots
\end{equation*}
The extra factor of $2$ comes from the fact that there are two values $\pm \sqrt{1 - |\zeta|^2}$ on the unit sphere for each $\zeta$ with $|\zeta| < 1$.
It is now clear that the contribution of the omitted term is $O(t^{-n/2-k+1/2})$.
Thus, we have shown that
\begin{equation}
\tr P E(t) = (\frac{n}{2} + k-1)\dots (\frac{n}{2})  \frac{4 \,  t^{-\frac{n}{2}-k}       }{(4\pi)^{n/2} \vol(S^{n-1})} \int\limits_{B^* Y} \sigma(P) \, d\mu_B \  + O(t^{-\frac{n}{2}-k+\frac1{2}})
\label{P-trace}\end{equation}
as $t \to 0$.

This result holds in particular for multiplication operators $\phi$, supported away from $\Sigma$:
\begin{equation}
\tr \phi E(t) \sim \frac{2 \,  t^{-\frac{n}{2}}}{(4\pi)^{n/2}} \int\limits_{ Y} \phi \, d\sigma  + O(t^{-\frac{n}{2}+\frac1{2}}).
\label{phi-tr}\end{equation}
The Gaussian bound \eqref{gaussian} shows that this holds for \emph{all} multiplication operators $\phi$, not necessarily supported away from the singular set. Hence \eqref{P-trace} holds also for $P$ equal to the identity operator:
\begin{equation}
\tr E(t) = \frac{2 t^{-n/2} \vol(Y)}{(4\pi)^{n/2}}  + o(t^{-\frac{n}{2}}), \ t \to 0.
\label{Id-trace}\end{equation}

\subsection{Dirichlet boundary condition}
The Seeley parametrix for the Dirichlet boundary condition on a manifold with smooth boundary still takes the form \eqref{Seel}, but the top terms $\gamma_0$ and $\delta_0$ now are
\begin{equation}\begin{aligned}
\gamma_0(x, y, \xi, \eta, t) &= e^{-t ( \xi^2 + | \eta |^2_h)} \\
\delta_0(x, y, \xi, \eta, t) &= - e^{-t (\xi^2 + | \eta |^2_h)} e^{-i\xi x}. \end{aligned}\label{top-terms-D}\end{equation}

Lemma~\ref{HK} gets replaced by
\begin{lem}\label{HK-D} Suppose that $\Omega \subset \RR^n$ is a
piecewise smooth manifold.
Then the local expansion of the Dirichlet heat kernel $e^{-tH_K}(z,z')$ is valid on any compact subset of the regular part $Y^o$ of the boundary.
\end{lem}

\begin{proof} This is proved as for Lemma~\ref{HK}, using the bound \begin{equation}
| e^{-t\Delta_{\Dir}}(z,z') | \leq (4\pi t)^{-n/2} e^{-|z-z'|^2/4t}.
\label{gaussian-D}\end{equation}
which follows by comparison with the heat kernel on $\RR^n$. 
\end{proof}

Next, we prove an identity for Dirichlet eigenfunctions on Lipschitz domains which is elementary for manifolds with smooth boundary.

\begin{lem} The Dirichlet eigenfunctions $u_j$ satisfy the equation
\begin{equation}
2 \lambda_j^2 = \int_Y (x \cdot \nu) |u_j^b|^2 \, d\sigma.
\label{Rellich}\end{equation}
\end{lem}

\begin{proof} For a smooth domain this formula was first proved by Rellich; see \cite{HT} for a proof. For a Lipschitz domain, the conclusion is still valid; see \cite{Verchota}, sections 1 and 2. 
\end{proof}

A similar proof, with the vector field $X$ replaced by one which is pointing out of the domain everywhere, shows that
$$
\| u_j^b \|^2 \leq C \lambda_j^2.
$$
It follows from this that the kernel $\tilde E(t)$ is trace class, the argument being the same as in the Neumann case.

Next we show

\begin{lem}\label{tEt} Let $P$ be either the identity operator on $Y$, or a differential operator of order $k$ on $Y$ supported away from $\Sigma$, with principal symbol $p(y, \eta)$. Then the trace of $P \tilde E(t)$ satisfies
\begin{equation}
\tr P \tilde E(t) = (\frac{n}{2} + k) (\frac{n}{2} + k - 1) \dots \frac{n}{2}  \frac{4 t^{-\frac{n}{2} - k -1}}{(4\pi)^{\frac{n}{2}} \vol(S^{n-1})} \int\limits_{B^* Y} p \, d\mu_{\Dir}   + o(t^{-\frac{n}{2} - k -1}).
\label{tt}\end{equation}
\end{lem}

\begin{proof}
We first consider the case where the coefficients of $P$ are supported away from $\Sigma$.
We compute as for the Neumann case, using the expression \eqref{Seel} for the heat kernel, with terms $\gamma_0$, $\delta_0$ now given by \eqref{top-terms-D}, and remembering to take normal derivatives in both variables. We obtain
$$
\tr P \tilde E(t) =  (2\pi)^{-n} \int e^{-t ( \eta^2 + |\eta |^2_h)}  \Big( \tilde p_{2k}(y, \zeta) + \tilde p_{2k-1}(y, \zeta) \Big) 2\xi^2 \, \mu\, dy \, d\xi \, d\eta. $$
Following the Neumann calculation, the $\tilde p_{2k}$ term gives \eqref{tt}, while the $\tilde p_{2k-1}$ term gives an $O(t^{-n/2-k-1/2})$ error. Note that the effect of the Dirichlet boundary condition is to change $k$ to $k+1$, and to introduce the factor $\xi^2$ into the integral above. This $\xi^2$ factor becomes $1 - |\zeta|^2$ in the integral over $|\zeta| \leq 1$, changing the Neumann measure into the Dirichlet measure.

The result for the identity operator is a special case of the following lemma.
\end{proof}

\begin{lem}\label{A-trace} For every measurable set $A \subset Y$,
$$
\lim_{\lambda \to \infty} \frac1{N(\lambda)} \lambda_j^{-2} \sum_{\lambda_j
< \lambda} \int_A  |u_j^{b}|^2 \, d\sigma = \frac{2 \vol(A)}{n \vol(\Omega)}.
$$
\end{lem}

\noindent{\it Remark. } This result was proved for domains with smooth boundary in \cite{Ozawa}.

\begin{proof}
By the Karamata Tauberian theorem, it is sufficient to show that
$$
\sum_j e^{-t \lambda_j^2} \lambda_j^{-2} \int_A  |u_j^{b}|^2 \, d\sigma \sim   \frac{2 \vol(A)}{(4\pi)^{n/2} n} t^{-n/2}, \quad t \to 0.
$$
By integrating in $t$ down from $t=1$, it is sufficient to prove that
$$
\sum_j e^{-t \lambda_j^2} \int_A  |u_j^{b}|^2 \, d\sigma \sim  \frac{\vol(A)}{(4\pi)^{n/2}} t^{-n/2-1}, \quad t \to 0.
$$
This is the same thing as asking whether
\begin{equation}
\int_A  \tilde E(t, y, y) \, d\sigma(y) \sim \frac{\vol(A)}{(4\pi)^{n/2}}  t^{-n/2-1}, \quad t \to 0.
\label{whether}\end{equation}
Let us consider the quantity
\begin{equation}
\int_A (x \cdot \nu) \tilde E(t, y, y) \, d\sigma(y).
\label{Aq}\end{equation}
By \eqref{tt}, if $A$ is compactly contained in $Y^o$, then \eqref{Aq} has asymptotics
\begin{equation}
  \frac{2n t^{-\frac{n}{2} - k -1}}{(4\pi)^{\frac{n}{2}} \vol(S^{n-1})} \int\limits_{B^* A} (x \cdot \nu) d\mu_{\Dir}   + o(t^{-\frac{n}{2} - k -1}).
\label{tttt}\end{equation}
Since
$$
\int_{0}^1 \sqrt{1 - r^2} \, r^{n-2}\, dr = \frac{1}{n} \int_{0}^1 \frac1{\sqrt{1 - r^2}} \, r^{n-2} \, dr,
$$
$\omega_{\Dir}(1) = \omega_{\Neu}(1)/n$, so \eqref{tttt} is given by
\begin{equation}
(4\pi)^{-\frac{n}{2}} t^{-\frac{n}{2} - k -1}  \int_A (x \cdot \nu) \, d\sigma
+ o(t^{-\frac{n}{2} - k -1}).
\end{equation}

Next consider the case when $A = Y$. Then using \eqref{Rellich}, \eqref{Aq}  is equal to
\begin{equation}
\sum_j e^{-t \lambda_j^2} \int_Y (x \cdot \nu)  |u_j^{b}|^2 \, d\sigma
= 2 \sum_j e^{-t \lambda_j^2} \lambda_j^2 \langle u_j, u_j \rangle
= 2 \tr \Delta e^{-t \Delta}.
\end{equation}
This has asymptotics
\begin{multline}
 2 \tr \Delta e^{-t \Delta} = -2 \frac{d}{dt} \tr e^{-t \Delta} = \frac{n \vol(\Omega) }{ (4\pi)^{n/2}}  \, t^{-n/2 - 1}  + o(t^{-\frac{n}{2} - k -1}) \\ = (4\pi)^{-\frac{n}{2}} \bigg( \int_Y x \cdot \nu \, d\sigma \bigg) t^{-n/2 - 1} + o(t^{-\frac{n}{2} - k -1}), \quad t \to 0.
\label{trr}\end{multline}

Now suppose for a moment that $\Omega$ is \emph{starshaped}. Then we can choose an origin so that $x \cdot \nu $ is positive on $Y$, and this together with the positivity of $\tilde E(t, y, y)$, the fact that $Y^o$ has full measure, and \eqref{trr} implies that \eqref{whether} holds for all measurable $A$.

Now we note that property \eqref{whether} depends only the asymptotics of $\tilde E(t, y, y)$ as $t \to 0$. This is local in $y$, so \eqref{whether} therefore holds for \emph{locally} starshaped domains. But this includes all Lipschitz domains, so we have finished the proof of the lemma.
\end{proof}

\subsection{Commutators of pseudodifferential operators} In this subsection we prove a Lemma needed in section~\ref{Psi}.

\begin{lem}\label{commuting}
Let $K$ be a homogeneous pseudodifferential operator of order $m$ and $Q_h$ a semiclassical pseudodifferential operator on $Y$ of order $0$, both supported away from $\Sigma$ in both variables. Assume also that $Q_h$ has compact microsupport. Then the operator norm of $[K, Q_h]$ is $o(h^{-m})$ as $h \to 0$.
\end{lem}

\begin{proof} We prove this in the case of operators on a manifold without boundary. The result above follows by localization to a neighbourhood of the support of $K$. Thus, we assume now that $\Sigma$ is empty. Also, we restrict to the case $m=0$, since the general case follows with minor modifications to the argument.

Let $\Delta$ be some positive elliptic differential operator of order $2$ on $Y$. We first prove the result when $Q_h = \phi(h^2 \Delta)$ for some smooth function $\phi$ with compact support. One can expect the theorem to be easier in this case since $Q_h$ is a semiclassical operator, but closely related to the homogeneous operator $\Delta$ which has well-behaved commutators with $K$.
We use the commutator formula and an almost analytic extension $\Phi$ of $\phi$ to write the commutator
\begin{equation*}\begin{gathered}
 {} [K, \phi(h^2 \Delta)] \\  = [K, h^2 \Delta] \phi'(h^2 \Delta) +
\frac{1}{2\pi} \int_{\mathbb{C}} \overline{\partial} \Phi(z) (h^2 \Delta - z)^{-1} [[K, h^2 \Delta], h^2 \Delta] (h^2 \Delta - z)^{-2} \, dz d\overline{z}.
\end{gathered}\end{equation*}
which follows from the Helffer-Sj\"ostrand formula \cite{HS}
$$
\phi(h^2 \Delta) = \frac{1}{2\pi} \int_{\mathbb{C}} \overline{\partial} \Phi(z) (h^2 \Delta - z)^{-1}  \, dz d\overline{z}.
$$
The first term may be written
\begin{equation}\begin{gathered}
h^2 [K, \Delta] (1 + h^2 \Delta)^{-1/2}  (1 + h^2 \Delta)^{1/2} \phi'(h^2 \Delta) \\
= h^2  [K, \Delta] \Delta^{-1/2}  \Delta^{1/2} (1 + h^2 \Delta)^{-1/2} \tilde \phi(h^2 \Delta) \\
= h \Big( [K, \Delta]  \Delta^{-1/2} \Big) \Big(  h\Delta^{1/2} (1 + h^2 \Delta)^{-1/2}\Big) \tilde \phi(h^2 \Delta).
\end{gathered}\label{firstterm}\end{equation}
The first factor in large parentheses is a pseudodifferential operator of order zero, the second is an operator with operator norm bounded by $1$ and the third is a semiclassical operator of order zero. Hence the operator norm of the \eqref{firstterm} is $O(h)$.

To deal with the integral over $\mathbb{C}$, we insert the factor
$\Delta \Delta^{-1} = (\Delta - h^{-2}z + h^{-2} z) \Delta^{-1}$, and get
$$\begin{gathered}
\frac{1}{\pi} \int_{\mathbb{C}} \overline{\partial} \Phi(z) (h^2 \Delta - z)^{-1} \Delta \Delta^{-1} [[K,  h^2 \Delta],  h^2 \Delta] (h^2 \Delta - z)^{-2} \, dz d\overline{z} \\
= Ch^2  \int_{\mathbb{C}} \overline{\partial} \Phi(z) \Big( \Id + z(h^2\Delta - z)^{-1} \Big) \Big( \Delta^{-1} [[K, \Delta],  \Delta] \Big) (h^2 \Delta - z)^{-2} \, dz d\overline{z}
\end{gathered}$$
Recall that $|\overline{\partial} \Phi(z)| \leq C_N |\operatorname{Im} z|^{-N-1} \langle z \rangle^{-N}$ for any $N$. Hence we can estimate the operator norm of this integral by
$$
Ch^2  \int_{\mathbb{C}} \overline{\partial} \Phi(z)
(1 + \langle z \rangle |\operatorname{Im} z|^{-1} ) |\operatorname{Im} z|^{-2}
\, dz d\overline{z}
\leq C h^2 \langle z \rangle^{-3} \leq C h^2.
$$
This proves the Lemma when $Q_h = \phi(h^2 \Delta)$.

To prove the Lemma in general, we choose a function $\phi$ so that $\phi(|\eta|^2)$ is identically $1$ on the microsupport of $Q_h$. This is possible since $Q_h$ has compact microsupport by hypothesis. Then, for any $\epsilon > 0$, we can find an $h_0 > 0$ and a differential operator $P$, or order $m$ say, such that the operator norm
$$
\| Q_h - h^m P \phi(h^2 \Delta) \|_{L^2(Y) \to L^2(Y)} \leq \epsilon \text{ for all } h \leq h_0.
$$
This follows from the density of polynomials in the $C^\infty$ topology on compact subsets, as in the proof of Lemma~\ref{LWL} in Section~\ref{Weyl}. Hence, it is enough to prove the result for all operators of the form $h^m P \phi(h^2 \Delta)$. (This is rather similar in spirit to the usual proof of the Riemann-Lebesgue lemma.)

To do this, we write
$$
[K, h^m P \phi(h^2 \Delta)] = h^m P [K, \phi(h^2 \Delta)] + [K, h^m P] \phi(h^2 \Delta),
$$
and consider each term. The second term is of the form
\begin{equation}
h^m A \, \phi(h^2 \Delta), \text{ where } A \text{ is a pseudo of order } m-1,
\label{psm}\end{equation}
 which can be treated as in \eqref{firstterm}.
The first we expand as
$$\begin{aligned}
&h^m P [K, \phi(h^2 \Delta)] \\
&= h^m P (1 + \Delta)^{-m/2} (1 + \Delta)^{m/2} (1 + h^2\Delta)^{-m/2} (1 + h^2\Delta)^{-m/2} [K, \phi(h^2 \Delta)] \\
&= \Big( P  (1 + \Delta)^{-m/2} \Big)\Big( (h^2 + h^2 \Delta)^{m/2}(1 + h^2 \Delta)^{-N} \Big) \Big(  (1 + h^2\Delta)^{N} [K, \phi(h^2 \Delta)] \Big).
\end{aligned}$$
The first two factors in large parentheses are bounded operators uniformly in $h$, provided $N > m/2$. For the remaining factor we write
$$
(1 + h^2\Delta)^{N} [K, \phi(h^2 \Delta)]  = [K, (1 + h^2\Delta)^{N} \phi(h^2 \Delta)] -  [K, (1 + h^2\Delta)^{N}] \phi(h^2 \Delta)
$$
and note the the first term has operator norm $O(h)$ by the first part of the proof, while the second is a sum of terms of the form \eqref{psm}. This completes the proof of the Lemma.
\end{proof}

\end{document}